\documentclass{elsarticle}
\makeatletter
\def\ps@pprintTitle{
  \let\@oddfoot\@empty
}
\makeatother

\usepackage[margin=1in]{geometry}

\usepackage{bbm}
\usepackage{graphicx}
\usepackage{pstool}
\usepackage{wrapfig}
\usepackage{caption}
\usepackage{subcaption}
\usepackage[section]{placeins} 
\usepackage{fancyhdr} 
\usepackage{lastpage} 
\usepackage{extramarks} 

\usepackage{xcolor} 
\usepackage{enumerate} 
\usepackage{paralist} 
\usepackage{amsmath, amsthm, amssymb, mathtools}
\usepackage{mathabx, pifont, stmaryrd} 
\usepackage[explicit]{titlesec} 
\usepackage{etoolbox} 
\usepackage{bibentry} 
\makeatletter\let\saved@bibitem\@bibitem\makeatother 
\usepackage[colorlinks, bookmarksopen, bookmarksnumbered,
            citecolor=red,urlcolor=red]{hyperref} 
\makeatletter\let\@bibitem\saved@bibitem\makeatother 
\usepackage{comment} 

\usepackage{algorithm}
\usepackage{algorithmic}

\usepackage{optidef}
\usepackage{mathtools}

\newtheorem{remark}{Remark}

\theoremstyle{definition}

\usepackage{bm} 
\usepackage{amsfonts} 

\newcommand{\Fcal}{\ensuremath{\mathcal{F}}}

\newcommand{\Hcal}{\ensuremath{\mathcal{H}}}

\newcommand{\Mcal}{\ensuremath{\mathcal{M}}}

\newcommand{\Nbb}{\ensuremath{\mathbb{N} }}

\newcommand{\Rbb}{\ensuremath{\mathbb{R} }}

\usepackage{tikz}
\usepackage{pgfplots}
\usepackage{xcolor}
\usepackage{pgfplotstable, booktabs}
        
\pgfplotsset{compat=1.9}
\definecolor{lightergreen}{rgb}{0.8, 1.0, 0.8}

\usetikzlibrary{pgfplots.groupplots}
\usepgfplotslibrary{fillbetween}
\usetikzlibrary{calc,fit,matrix,arrows,automata,positioning,shapes}
\usetikzlibrary{arrows.meta}

\pgfplotsset{select coords between index/.style 2 args={
    x filter/.code={
        \ifnum\coordindex<#1\fi
        \ifnum\coordindex>#2\fi
    }
}}

\tikzset{
 invisible/.style={opacity=0},
 visible on/.style={alt={#1{}{invisible}}},
 alt/.code args={<#1>#2#3}{%
   \alt<#1>{\pgfkeysalso{#2}}{\pgfkeysalso{#3}}
 },
}

\newcommand{\colorbarMatplotlibCoolwarm}[7]{
\begin{tikzpicture}
\begin{axis}[
   hide axis, scale only axis,
   xmin=0, xmax=1, 
   ymin=0, ymax=1, 
   height=0pt, width=0pt,
   colormap={coolwarm}{rgb255=(59,76,192) rgb255=(60,78,194) rgb255=(61,80,195) rgb255=(62,81,197) rgb255=(63,83,198) rgb255=(64,85,200) rgb255=(66,87,201) rgb255=(67,88,203) rgb255=(68,90,204) rgb255=(69,92,206) rgb255=(70,93,207) rgb255=(71,95,209) rgb255=(73,97,210) rgb255=(74,99,211) rgb255=(75,100,213) rgb255=(76,102,214) rgb255=(77,104,215) rgb255=(79,105,217) rgb255=(80,107,218) rgb255=(81,109,219) rgb255=(82,110,221) rgb255=(84,112,222) rgb255=(85,114,223) rgb255=(86,115,224) rgb255=(87,117,225) rgb255=(89,119,226) rgb255=(90,120,228) rgb255=(91,122,229) rgb255=(93,123,230) rgb255=(94,125,231) rgb255=(95,127,232) rgb255=(96,128,233) rgb255=(98,130,234) rgb255=(99,131,235) rgb255=(100,133,236) rgb255=(102,135,237) rgb255=(103,136,238) rgb255=(104,138,239) rgb255=(106,139,239) rgb255=(107,141,240) rgb255=(108,142,241) rgb255=(110,144,242) rgb255=(111,145,243) rgb255=(112,147,243) rgb255=(114,148,244) rgb255=(115,150,245) rgb255=(116,151,246) rgb255=(118,153,246) rgb255=(119,154,247) rgb255=(120,156,247) rgb255=(122,157,248) rgb255=(123,158,249) rgb255=(124,160,249) rgb255=(126,161,250) rgb255=(127,163,250) rgb255=(129,164,251) rgb255=(130,165,251) rgb255=(131,167,252) rgb255=(133,168,252) rgb255=(134,169,252) rgb255=(135,171,253) rgb255=(137,172,253) rgb255=(138,173,253) rgb255=(140,174,254) rgb255=(141,176,254) rgb255=(142,177,254) rgb255=(144,178,254) rgb255=(145,179,254) rgb255=(147,181,255) rgb255=(148,182,255) rgb255=(149,183,255) rgb255=(151,184,255) rgb255=(152,185,255) rgb255=(153,186,255) rgb255=(155,187,255) rgb255=(156,188,255) rgb255=(158,190,255) rgb255=(159,191,255) rgb255=(160,192,255) rgb255=(162,193,255) rgb255=(163,194,255) rgb255=(164,195,254) rgb255=(166,196,254) rgb255=(167,197,254) rgb255=(168,198,254) rgb255=(170,199,253) rgb255=(171,199,253) rgb255=(172,200,253) rgb255=(174,201,253) rgb255=(175,202,252) rgb255=(176,203,252) rgb255=(178,204,251) rgb255=(179,205,251) rgb255=(180,205,251) rgb255=(182,206,250) rgb255=(183,207,250) rgb255=(184,208,249) rgb255=(185,208,248) rgb255=(187,209,248) rgb255=(188,210,247) rgb255=(189,210,247) rgb255=(190,211,246) rgb255=(192,212,245) rgb255=(193,212,245) rgb255=(194,213,244) rgb255=(195,213,243) rgb255=(197,214,243) rgb255=(198,214,242) rgb255=(199,215,241) rgb255=(200,215,240) rgb255=(201,216,239) rgb255=(203,216,238) rgb255=(204,217,238) rgb255=(205,217,237) rgb255=(206,217,236) rgb255=(207,218,235) rgb255=(208,218,234) rgb255=(209,219,233) rgb255=(210,219,232) rgb255=(211,219,231) rgb255=(213,219,230) rgb255=(214,220,229) rgb255=(215,220,228) rgb255=(216,220,227) rgb255=(217,220,225) rgb255=(218,220,224) rgb255=(219,220,223) rgb255=(220,221,222) rgb255=(221,221,221) rgb255=(222,220,219) rgb255=(223,220,218) rgb255=(224,219,216) rgb255=(225,219,215) rgb255=(226,218,214) rgb255=(227,218,212) rgb255=(228,217,211) rgb255=(229,216,209) rgb255=(230,216,208) rgb255=(231,215,206) rgb255=(232,215,205) rgb255=(232,214,203) rgb255=(233,213,202) rgb255=(234,212,200) rgb255=(235,212,199) rgb255=(236,211,197) rgb255=(236,210,196) rgb255=(237,209,194) rgb255=(238,209,193) rgb255=(238,208,191) rgb255=(239,207,190) rgb255=(240,206,188) rgb255=(240,205,187) rgb255=(241,204,185) rgb255=(241,203,184) rgb255=(242,202,182) rgb255=(242,201,181) rgb255=(243,200,179) rgb255=(243,199,178) rgb255=(244,198,176) rgb255=(244,197,174) rgb255=(245,196,173) rgb255=(245,195,171) rgb255=(245,194,170) rgb255=(245,193,168) rgb255=(246,192,167) rgb255=(246,191,165) rgb255=(246,190,163) rgb255=(246,188,162) rgb255=(247,187,160) rgb255=(247,186,159) rgb255=(247,185,157) rgb255=(247,184,156) rgb255=(247,182,154) rgb255=(247,181,152) rgb255=(247,180,151) rgb255=(247,178,149) rgb255=(247,177,148) rgb255=(247,176,146) rgb255=(247,174,145) rgb255=(247,173,143) rgb255=(247,172,141) rgb255=(247,170,140) rgb255=(247,169,138) rgb255=(247,167,137) rgb255=(247,166,135) rgb255=(246,164,134) rgb255=(246,163,132) rgb255=(246,161,131) rgb255=(246,160,129) rgb255=(245,158,127) rgb255=(245,157,126) rgb255=(245,155,124) rgb255=(244,154,123) rgb255=(244,152,121) rgb255=(244,151,120) rgb255=(243,149,118) rgb255=(243,147,117) rgb255=(242,146,115) rgb255=(242,144,114) rgb255=(241,142,112) rgb255=(241,141,111) rgb255=(240,139,109) rgb255=(240,137,108) rgb255=(239,136,106) rgb255=(238,134,105) rgb255=(238,132,103) rgb255=(237,130,102) rgb255=(236,129,100) rgb255=(236,127,99) rgb255=(235,125,97) rgb255=(234,123,96) rgb255=(233,121,95) rgb255=(233,120,93) rgb255=(232,118,92) rgb255=(231,116,90) rgb255=(230,114,89) rgb255=(229,112,88) rgb255=(228,110,86) rgb255=(227,108,85) rgb255=(227,106,83) rgb255=(226,104,82) rgb255=(225,102,81) rgb255=(224,100,79) rgb255=(223,98,78) rgb255=(222,96,77) rgb255=(221,94,75) rgb255=(220,92,74) rgb255=(218,90,73) rgb255=(217,88,71) rgb255=(216,86,70) rgb255=(215,84,69) rgb255=(214,82,67) rgb255=(213,80,66) rgb255=(212,78,65) rgb255=(210,75,64) rgb255=(209,73,62) rgb255=(208,71,61) rgb255=(207,69,60) rgb255=(205,66,59) rgb255=(204,64,57) rgb255=(203,62,56) rgb255=(202,59,55) rgb255=(200,57,54) rgb255=(199,54,53) rgb255=(198,51,52) rgb255=(196,49,50) rgb255=(195,46,49) rgb255=(193,43,48) rgb255=(192,40,47) rgb255=(190,37,46) rgb255=(189,34,45) rgb255=(188,30,44) rgb255=(186,26,43) rgb255=(185,22,41) rgb255=(183,17,40) rgb255=(181,11,39) rgb255=(180,4,38)},
   colorbar horizontal,
   point meta min=#1, point meta max=#5,
   colorbar style={width=#6, height=#7,xtick={#1,#2,#3,#4,#5},ticklabel style={font=\normalsize}}
]
\addplot [draw=none] coordinates {(0,0)};
\end{axis}\textbf{}
\end{tikzpicture}
}

\usepackage{setspace}
\usepackage{enumitem}
\usepackage{ulem}
\usepackage{tipa}

\newbool{fastcompile}
\setbool{fastcompile}{false}

\begin{document}
\title{Neural network-based Godunov corrections for approximate Riemann solvers using bi-fidelity learning}

\author[rvt2]{Akshay Thakur\fnref{fn1}}
\ead{athakur3@nd.edu}

\author[rvt2]{Matthew J. Zahr\fnref{fn2}\corref{cor1}}
\ead{mzahr@nd.edu}

\address[rvt2]{Department of Aerospace and Mechanical Engineering, University
	of Notre Dame, Notre Dame, IN 46556, United States}

\cortext[cor1]{Corresponding author}

\fntext[fn1]{Graduate Student, Department of Aerospace and Mechanical Engineering, University of Notre Dame}
\fntext[fn2]{Assistant Professor, Department of Aerospace and Mechanical
             Engineering, University of Notre Dame}


\begin{abstract}
The Riemann problem is fundamental in the computational modeling of hyperbolic partial differential equations, enabling the development of stable and accurate upwind schemes. While exact solvers provide robust upwinding fluxes, their high computational cost necessitates approximate solvers. Although approximate solvers achieve accuracy in many scenarios, they produce inaccurate solutions in certain cases. To overcome this limitation, we propose constructing neural network-based surrogate models, trained using supervised learning, designed to map interior and exterior conservative state variables to the corresponding exact flux. Specifically, we propose two distinct approaches: one utilizing a vanilla neural network and the other employing a bi-fidelity neural network. The performance of the proposed approaches is demonstrated through applications to one-dimensional and two-dimensional partial differential equations, showcasing their robustness and accuracy.
\end{abstract}
    
\maketitle

\section{Introduction}
\label{sec:intro}
In computational modeling of systems governed by hyperbolic partial differential equations (PDEs), the Riemann problem and its solution play a pivotal role, particularly in the development of upwind schemes \cite{godunov1959finite,toro2013riemann}, which ensure stability, provided the Courant-Friedrichs-Lewy condition is satisfied, and accuracy in computing the numerical solution for conservation laws. The Riemann problem can be solved using an exact (Godunov) or an approximate solver, such as Roe \cite{roe1981approximate} and  Harten-Lax-van Leer (HLL) \cite{harten1983upstream} solvers. Approximate Riemann solvers are widely favored in numerical methods due to their computational efficiency as compared to the exact solvers. However, they come with some drawbacks. Specifically, fluxes obtained from approximate solvers are generally less accurate and more diffusive than an exact Riemann solver, leading to the loss or smoothing of fine wave structures and sharp features \cite{ketcheson2020riemann}. This issue is aggravated further when working with coarser meshes or low-order discretizations. Additionally, approximate solvers such as Roe and HLL struggle with strong shocks, transonic rarefactions, and complex wave interactions \cite{toro2013riemann}. Although the difficulties that the Roe solver faces with transonic rarefactions can be remedied using an entropy fix, it fails to preserve positivity (with or without entropy fix) and faces robustness issues with near dry states or low density flows \cite{ketcheson2020riemann}. Similarly, modifications can be made to the HLL solver to handle contact discontinuities, which again may enhance accuracy selectively but lacks universality.\par 
Recent years have witnessed a notable surge in the development and application of deep learning techniques for the construction of surrogate models (SMs) and augmentation of PDE solvers with neural networks (NNs) in complex physical and engineering applications. This includes closure modeling \cite{nair2023deep,sirignano2023deep,guan2022stable}, multi-fidelity learning \cite{meng2020composite}, physics-informed neural networks \cite{gao2022physics, raissi2019physics, wang2023fluxnet}, differentiable solvers \cite{nair2023deep, bezgin2023jax}, neural operators \cite{lu2021learning,li2020fourier}, and supervised neural network regression for specific modeling tasks \cite{bar2019learning, gyrya2024machine,guan2022stable}. Concurrently, the exploration of the potential of neural networks to augment or replace Riemann solvers, to provide effective approximations while preserving critical attributes like resilience in the face of strong shocks or transonic rarefactions, has also garnered attention \cite{wang2023fluxnet, bezgin2023jax, gyrya2024machine,magiera2020constraint,xu2024unsupervised,tong2024roenet}. Magiera et al. \cite{magiera2020constraint} trained a fully connected network with a penalty based on Rankine-Hugonoit shock jump conditions to predict the state variables in the star-state region for Riemann problems with ideal thermodynamics. Wang and Hickey \cite{wang2023fluxnet} followed a similar approach for Riemann problems with non-ideal thermodynamics. Gyrya et al. \cite{gyrya2024machine} applied supervised learning techniques, including Gaussian Process (GP) and NN-based SMs, to predict star-state variables and various wave speeds for Riemann problems. They also employed these models for classifying wave structures arising from different initial Riemann data. While their work demonstrated the feasibility of data-driven surrogates for exact Riemann solvers, they did not integrate the trained models into a finite volume (FV) framework or incorporate physics-based constraints or priors into the learning process. Furthermore, all three studies \cite{magiera2020constraint,wang2023fluxnet,gyrya2024machine} focused on creating SMs from primitive interior and exterior states to the primitive star states. Bezgin et al. \cite{bezgin2023jax} replaced the scalar numerical viscosity in the Rusanov flux function with a neural network and performed embedded training with a coarse-grained FV-solver in a differentiable way. Tong et al. \cite{tong2024roenet} replace the flux Jacobian matrix with a matrix constructed using eigenvector and eigenvalue matrices learned using neural networks. Again, embedded training is performed within a first-order finite volume setting. However, training by embedding the flux framework in a specific discretization becomes dependent on it, and therefore, the learned flux functions may not generalize to any other discretization. Xu et al. \cite{xu2024unsupervised} predict the interfacial pressure in a multi-material Riemann problem using an NN-based unsupervised learning approach. We propose a different approach, where we construct an NN-based surrogate using supervised learning as an approximate map from the interior and exterior conservative variables to the corresponding Godunov flux. However, a vanilla NN trained only on collected data might not generalize well to unseen scenarios. Therefore, we also propose employing bi-fidelity learning approaches for NN-based corrections to fluxes obtained from approximate solvers. \par 
The remainder of this paper is structured as follows. Section \ref{sec:probstat} outlines the problem setting by introducing the governing equations and their discretization. It also discusses the Riemann problem, numerical flux formulation, and the use of surrogates as numerical flux functions. Section \ref{sec:prop} outlines the proposed surrogate models and the associated optimization framework. Section \ref{sec:num_exp} assesses the models through numerical examples including one-dimensional and two-dimensional hyperbolic and parabolic PDEs. Finally, Section \ref{sec:conclude} presents the concluding remarks.

\section{Problem statement}
\label{sec:probstat}
In this section, we first present the governing equations and semidiscretization for a general time-dependent system of conservation laws, which will serve as the basis for modeling unsteady, nonlinear phenomena (Section \ref{ssec:goveqn}). We then introduce the concept of the Riemann problem and its role in defining the numerical flux at inter-element interfaces in the discretized domain (Section \ref{ssec:rpnf}), and two commonly used approximate Riemann solver (Section \ref{ssec:arp}). Lastly, we outline the motivation behind developing neural network-based surrogate models (SM) to approximate the inviscid numerical flux (Section \ref{ssec:motps}).

\subsection{Governing equations} \label{ssec:goveqn}

Consider a general form for a system of $m$ conservation laws defined in a spatial domain $\Omega \subset \mathbb{R}^d$ over the time interval $\mathcal{T} = (0,T]$, which can be expressed as follows:
\begin{equation}
   \frac{\partial U}{\partial t} + \nabla \cdot F(U, \nabla U) = S(U, \nabla U), \quad U(\cdot, 0) = U_{0}(\cdot),
   \label{eqn:genconslaw}
\end{equation}
where $U \in \mathbb{R}^{m}$ is the vector of conserved variables implicitly defined as the solution of  (\ref{eqn:genconslaw}), $F : \mathbb{R}^m \times \mathbb{R}^{m \times d}  \rightarrow \mathbb{R}^{m\times d} $ is the flux function, and $S: \mathbb{R}^m \times \mathbb{R}^{m \times d}  \rightarrow \mathbb{R}^{m}$ represents the source term. For any $W \in  \mathbb{R}^m$ and  $\nabla W \in  \mathbb{R}^{m \times d}$, the flux function can be split into viscous and inviscid terms as $F(W, \nabla W) = F_{I}(W) - F_{V}(W,\nabla W)$, where $F_{I}: \mathbb{R}^m \rightarrow \mathbb{R}^{m\times d}$ is the inviscid flux function and $F_{V}: \mathbb{R}^m \times \mathbb{R}^{m \times d}  \rightarrow \mathbb{R}^{m\times d}$ is the viscous flux function. Therefore, (\ref{eqn:genconslaw}) can also be expressed as:
\begin{equation}
\frac{\partial U}{\partial t} + \nabla \cdot F_{I}(U) = \nabla \cdot F_{V}(U, \nabla U) + S(U, \nabla U), \quad U(\cdot, 0) = U_{0}(\cdot).
\label{eqn:genconslaw2}
\end{equation}

\noindent Furthermore, the Jacobian of the projected inviscid flux can be expressed as:
\begin{equation}
B : \mathbb{R}^m \times \mathbb{N}_{d} \to \mathbb{R}^{m \times m}, \quad B(W, N) = \frac{\partial [F_I(W) N]}{\partial W},
\end{equation}
where $\mathbb{N}_{d} :=  \{N \in \mathbb{R}^{d} \mid ||N|| = 1\}$. In addition, consider a discretization operator $\mathcal{M}$ and a method of lines-based semi-discretization $\xi_d$ on the domain $\Omega$ for the PDE in (\ref{eqn:genconslaw}). For cell-centered finite volume or discontinuous finite element methods, $\xi_d$ is a decomposition of the domain $\Omega$ into a collection of elements. Then, for a single element $K \in \xi_d$, (\ref{eqn:genconslaw}) can be written as:
\begin{equation}
   \frac{\partial U_K}{\partial t}  = \mathcal{M}_{K}(U_{\bar{K}}),
   \label{eqn:semidis}
\end{equation}
where  $U_{\bar{K}}$ are collections of state vectors from the elements $\bar{K}\subset\xi_d$ surrounding $K$ (inclusive) required to compute the right-hand side. Furthermore, $\mathcal{M}_K$ depends on the flux function $F$, source term $S$, inviscid numerical flux function $\mathcal{H}$, and the semi-discretization $\xi_d$.
In this work, we use a cell-centered finite volume method, in which case $\Mcal_K$ is defined as
\begin{equation}
    \Mcal_K(U_{\bar{K}}) =
    -\frac{1}{|K|}\sum_{j=1}^{f'} |\partial K_j| \left(\Hcal(U_{\partial K_j}^+, U_{\partial K_j}^-, N_{\partial K_j}) - F_V(U_{\partial K_j},Q_{\partial K_j})N_{\partial K_j}\right) + S(U_K, Q_{K}),
     \label{eqn:semidis1}
\end{equation}
where $U_K$ is the cell-average of $U$ over the element $K\in\xi_d$, $Q_K$ is an approximated reconstruction of $\nabla U$ at the centroid of $K$, $\partial K_j$ is the $j$th face of element $K$ (assuming each $K\in\xi_d$ is a polyhedra with $f'$ planar faces) with outward unit normal $N_{\partial K_j}\in\Nbb^d$, $U_{\partial K_j}^+$ ($U_{\partial K_j}^-$) is a reconstruction of $U$ at the centroid of the face $\partial K_j$ on the interior (exterior) of element $K$, $U_{\partial K_j}$ and $Q_{\partial K_j}$ are single-valued reconstructions at the centroid of the face $\partial K_j$ that approximate $U$ and $\nabla U$, respectively. All reconstructions use the cell-averages in the elements of $\bar{K}\subset\xi_d$. The inviscid numerical flux function, $\Hcal : \Rbb^m \times \Rbb^m \times \Nbb^d \rightarrow \Rbb^m$ with $\Hcal : (U^+, U^-, N) \mapsto \Hcal(U^+,U^-,N)$, is a single-valued approximation of $F_I(U) N$. Multidimensional numerical fluxes are usually constructed from rotational invariance of the inviscid flux function, i.e.,
        \begin{equation}
            F_I(U) \cdot N = T_s^{-1} F_I^x(T_s U),
            \label{eqn:fvm3}
        \end{equation}
where \( F_I^x \) is the augmented inviscid flux function in the \( x \)-direction, which includes the standard one-dimensional inviscid system for the given PDE, augmented with advection equations for the tangential velocity, and \( T_s \) is the rotation matrix that depends on the normal $N$. For further details, interested readers are referred to \cite{toro2013riemann}. 

\subsection{Riemann problem and Godunov solver}\label{ssec:rpnf}
A Riemann problem is an initial value problem for a hyperbolic PDE, characterized by an initial condition that consists of two distinct constant states separated by a discontinuity. The solution to this problem involves the propagation of waves, such as shock waves, contact discontinuities, rarefaction waves, or a combination of them, which emanate from the initial discontinuity and evolve over time. These waves can propagate in multiple directions, depending on the nature of the underlying PDE and the problem settings. Furthermore, the solution of the Riemann problem is used in (\ref{eqn:semidis1}) to define the inviscid numerical flux function at the inter-element boundaries. In addition, solutions to both exact or approximate Riemann problems can be used in practice to define the inviscid numerical flux function, depending on the desired trade-off between computational cost, accuracy, and robustness.\par
In a multidimensional setting, for a specific element \( K \) with a unit outward normal vector \( N \), we consider the Riemann problem associated with the inviscid flux projected along the normal direction, denoted as \( F_{I} \cdot N \), given by

\begin{equation}
 \frac{\partial \Psi}{\partial \tau} + \frac{\left[ F_{I}(\Psi) \cdot N \right]}{\partial s} = 0, \quad \Psi(s,0) = 
\begin{cases}
U^+ & \text{if } s < 0, \\
U^- & \text{if } s \geq 0,
\end{cases}
\label{eqn:1drp}
\end{equation}
where \( \Psi(s, \tau)\) represents the solution of the Riemann problem, and \( U^+ \) and \( U^- \) are the interior and exterior traces of \( U \) at a point \( x \in \partial K \) and time \( t \). Due to the self-similarity of the Riemann problem solution, \( \Psi(s, \tau)\) can be expressed as a function of similarity variable \(s/\tau\), denoted as \(\tilde{\Psi}(s/\tau) \). Using the same property, we define the Godunov (exact) inviscid numerical flux by evaluating \(\tilde{\Psi}(s/\tau) \) at $s/\tau = 0$  as
\begin{equation}
\mathcal{H}_{gdnv}(U^+, U^-, N) = F_I(\tilde{\Psi}(0;N))\cdot N,
\label{eqn:exact}
\end{equation}
\noindent where the dependence of $\tilde{\Psi}$ on the unit normal $N$ indicates that the Riemann solution was obtained along this specific direction.
\subsection{Approximate solvers}\label{ssec:arp}
Although the inviscid numerical flux at an element interface can be determined by solving a one-dimensional Riemann problem exactly, this approach can become computationally prohibitive for complex nonlinear systems, such as the shallow water or Euler equations, particularly for large-scale simulations. To mitigate this challenge, approximate Riemann solvers were developed \cite{roe1981approximate, harten1983upstream}. These solvers provide an efficient alternative by simplifying the Riemann problem so that it can be solved exactly with a much lower computational effort. Among the various approximate solvers, one of the most widely used is the Roe solver \cite{roe1981approximate}. The Roe solver approximates the Riemann problem by linearizing the nonlinear flux function around an averaged state. Roe's numerical flux function is
\begin{equation}
\mathcal{H}_{roe}(U^+, U^-, N) = \frac{1}{2} \left[ F_I(U^+)\cdot N + F_I(U^-) \cdot N\right] - \frac{1}{2}|\tilde{B}(U^+, U^-, N)| (U^- - U^+),
\label{eqn:roe1}
\end{equation}

\noindent where \( \tilde{B}: \mathbb{R}^m \times \mathbb{R} \times \mathbb{N}^d \to \mathbb{R}^m\) is the projected inviscid flux Jacobian evaluated at Roe averaged states, i.e., 

\begin{equation}
\tilde{B}(U^+, U^-, N) = B(\hat{W}(U^+, U^-), N).
\label{eqn:roe2}
\end{equation}

\noindent \( \hat{W}: \mathbb{R}^m \times \mathbb{R}^m \to \mathbb{R}^m \) is the Roe average, which depends on the specific equations being solved and is chosen to make $\mathcal{H}_{roe}$ conservative at shocks. The absolute value of the square matrix \( |A| \) is defined through its eigenvalue decomposition as

\begin{equation}
|A| = L |\Lambda| L^{-1},
\end{equation}

\noindent where \( L \) and \( L^{-1} \) are the matrices of right and left eigenvectors of \( A \), respectively, and \( |\Lambda| \) is a diagonal matrix containing the absolute values of the eigenvalues of \( A \).\par Another popular approximate Riemann solver is the HLL solver.
The HLL solver is a two-wave solver that permits the left-moving and right-moving waves to propagate at different speeds. The HLL numerical flux is defined as
\begin{equation}
    \mathcal{H}_{hll}(U^+,U^-,N) =
    \begin{cases}
    \displaystyle F_I(U^+)\cdot N & \text{if } 0 \le S_+, \\
    \displaystyle \frac{S_- F_I(U^+)\cdot N - S_+ F_I(U^-)\cdot N + S_+ S_- (U^- - U^+)}{S_- - S_+} & \text{if } S_+ < 0 < S_-, \\
    \displaystyle F_I(U^-)\cdot N & \text{if } S_- \le 0,
    \end{cases}
    \label{eqn:hll1}
\end{equation}
where $S_-$ and $S_+$ are the two wave speeds, for which different choices have been proposed in the past literature. In this work, we compute the wave speeds based on the suggestion of Einfeldt \cite{einfeldt1988godunov}.

\subsection{Neural networks as an alternative} \label{ssec:motps}
The Godunov solver undoubtedly offers a more accurate and robust treatment of nonlinear systems, particularly for complex wave interactions, as compared to approximate solvers. However, these advantages come with a higher computational cost.
One promising direction that could bridge the gap between approximate and exact Riemann solvers is the development of NN-based surrogates for the Godunov solver. This approach involves constructing a neural network-based nonlinear map from the interior and exterior states $U^+$ and $U^-$, along with the normal $N$, to approximate the flux $\mathcal{H}(U^+, U^-, N)$ produced by the Godunov solver. Mathematically, this SM can be represented as

\begin{equation}
    \hat{\mathcal{H}}(U^+, U^-, N) = \mathcal{F}(U^+, U^-, N; \theta).
    \label{eqn:nnsm1}
\end{equation}

\noindent Here, $\mathcal{F}$ is a nonlinear approximator parameterized by $\theta$, and $\hat{\mathcal{H}}$ is the approximated Godunov flux function. In this study, our goal is to develop such surrogate approaches. Specifically, we aim to incorporate physics-based inductive biases from approximate solvers into NN-based models to create robust and accurate surrogates for the Godunov flux.


\section{Proposed models}
\label{sec:prop}
In this section, we present two distinct models that will be utilized as surrogates for the Godunov flux: a vanilla neural network (Section~\ref{ssec:models:vanilla}) and a neural network embedded in a bi-fidelity learning strategy (Section~\ref{ssec:models:bifidelity}). Both Godunov flux SMs are constructed in one dimension and generalized to higher dimensions using rotational invariance (\ref{eqn:fvm3}) as
\begin{equation}
 \Hcal(U^+, U^-, N) = T_s^{-1} \Hcal^x(T_s U^+, T_s U^-),
\end{equation}
where $\Hcal^x : \Rbb^m \times \Rbb^m\rightarrow\Rbb^m$ is a one-dimensional, $x$-directed, single-valued, numerical flux function. Finally, we detail our approach to train the neural networks used in the two models (Section~\ref{ssec:ohd}).

\subsection{A vanilla neural network as a Godunov flux surrogate} \label{ssec:models:vanilla}
As our first model, we utilize a fully connected neural network-based (FCNN-based) surrogate to approximate the one-dimensional flux obtained from the Godunov solver.  Mathematically, this can be represented by modifying (\ref{eqn:nnsm1}) as follows:

\begin{equation}
    \mathcal{\hat{H}}_{vn}^x(U^+, U^-; \theta) = \mathcal{F}_{NN}(U^+, U^-; \theta),
    \label{eqn:nnsm2}
\end{equation}

\noindent where \( \mathcal{F}_{NN} \) represents the neural network with parameters \( \theta \). Because we are using the $x$-directed numerical flux function $\Hcal_x$, the unit normal \( N \) does not appear. As such it will not be used to define or train the neural network $\Fcal_{NN}$.

\begin{figure}[ht!]
    \centering
    \begin{tikzpicture}[x=1.2cm, y=1.5cm, >=stealth]
    
        \node[draw, minimum size=0.6cm] (I-1) at (0, -1) {$U^+$};
        \node[draw, minimum size=0.6cm] (I-2) at (0, -2) {$U^-$};
    
        \foreach \i in {1,2,3}
            \node[draw, circle, fill=lightergreen, minimum size=0.6cm] (H1-\i) at (1.5, -\i + 0.5) {};
    
        \foreach \i in {1,2,3}
            \node[draw, circle, fill=lightergreen, minimum size=0.6cm] (H2-\i) at (3, -\i + 0.5) {};
    
        \foreach \i in {1,2,3}
            \node[draw, circle, fill=lightergreen, minimum size=0.6cm] (H3-\i) at (4.5, -\i + 0.5) {};
    
        \node[draw, minimum size=0.6cm] (O-1) at (6.0, -1.5) {$f^{H}$}; 
    
        \foreach \i in {1,2}
            \foreach \j in {1,2,3}
                \draw[->] (I-\i) -- (H1-\j);
    
        \foreach \i in {1,2,3}
            \foreach \j in {1,2,3}
                \draw[->] (H1-\i) -- (H2-\j);
    
        \foreach \i in {1,2,3}
            \foreach \j in {1,2,3}
                \draw[->] (H2-\i) -- (H3-\j);
    
        \foreach \i in {1,2,3}
            \draw[->] (H3-\i) -- (O-1);
    
    \end{tikzpicture}
    \caption{Schematic illustration of a fully connected neural network used as a surrogate for the Godunov flux.}
    \label{fig:nnvanilla}
\end{figure}
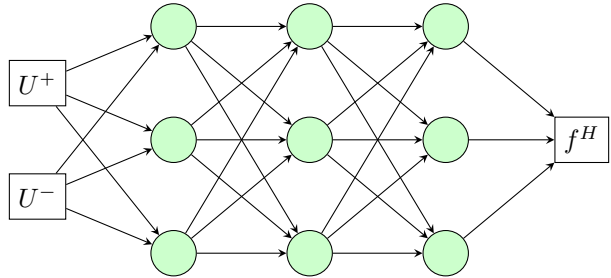

Figure~\ref{fig:nnvanilla} provides a schematic illustration of the vanilla NN (VNN) used as a surrogate for the Godunov flux. The network consists of input layers, hidden layers with nonlinear activations, and an output layer that produces the flux approximation. The inputs to the network are the interior $U^+$ and exterior $U^-$ states.
The NN-based surrogate is trained to approximate the flux using supervised learning, necessitating the generation of a suitable dataset. This dataset, denoted as \( \mathcal{D}_H \), is constructed using training data obtained from the Godunov solver as follows:

\begin{equation}
   \mathcal{D}_H = \left\{ \left( U^+_i, U^-_i \right), f_i^{H} = \mathcal{H}_{gdnv}(U^+_i, U^-_i,e_1) \right\}_{i=1}^{N_H},
   \label{eqn:data}
\end{equation}

\noindent where $e_1\in\Rbb^d$ is first canonical unit vector (we assume $\Hcal^x$ is directed in the positive $x$ direction) and \( f^H \) represents the flux values computed by the Godunov solver. This dataset comprises \( N_H \) input-output pairs, where each pair includes the inputs \( U^+ \) and \( U^- \), sampled from some suitable distribution as described in Section \ref{sec:num_exp}, and the corresponding Godunov flux \( f^H \). The goal of the NN surrogate is to learn a nonlinear approximation map from the input states in the dataset to the corresponding Godunov flux.

    \subsection{Bi-fidelity neural network as a Godunov flux surrogate} \label{ssec:models:bifidelity}
    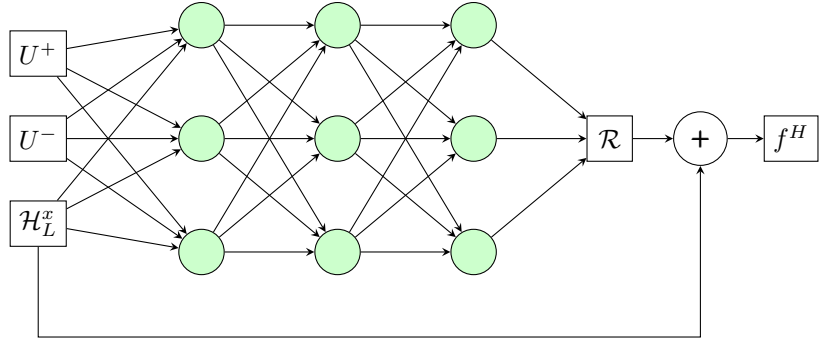
\begin{figure}[ht!]
        \centering
        \begin{tikzpicture}[x=1.2cm, y=1.5cm, >=stealth]
    
            \node[draw, minimum size=0.6cm] (I-1) at (-0.3, -0.75 ) {$U^+$};
            \node[draw, minimum size=0.6cm] (I-2) at (-0.3, -2*0.75 ) {$U^-$};
            \node[draw, minimum size=0.6cm] (I-3) at (-0.3, -3*0.75 ) {$\mathcal{H}_L^x$};
    
            \foreach \i in {1,2,3}
                \node[draw, circle, fill=lightergreen, minimum size=0.6cm] (H1-\i) at (1.5, -\i + 0.5) {};
    
            \foreach \i in {1,2,3}
                \node[draw, circle, fill=lightergreen, minimum size=0.6cm] (H2-\i) at (3, -\i + 0.5) {};
    
            \foreach \i in {1,2,3}
                \node[draw, circle, fill=lightergreen, minimum size=0.6cm] (H3-\i) at (4.5, -\i + 0.5) {};
    
            \node[draw, minimum size=0.6cm] (O-1) at (6.0, -1.5) {$\mathcal{R}$};
    
            \node[draw, circle, minimum size=0.6cm] (plus-circle) at (7.0, -1.5) {$\boldsymbol{+}$};
    
            \node[draw, minimum size=0.6cm] (final-output) at (8.0, -1.5) {$f^{H}$};
    
            \foreach \i in {1,2,3}
                \foreach \j in {1,2,3}
                    \draw[->] (I-\i) -- (H1-\j);
    
            \foreach \i in {1,2,3}
                \foreach \j in {1,2,3}
                    \draw[->] (H1-\i) -- (H2-\j);
    
            \foreach \i in {1,2,3}
                \foreach \j in {1,2,3}
                    \draw[->] (H2-\i) -- (H3-\j);
    
            \foreach \i in {1,2,3}
                \draw[->] (H3-\i) -- (O-1);
    
            \draw[->] (O-1) -- (plus-circle);
    
            \draw[->] (plus-circle) -- (final-output);
    
            \draw[->] (I-3) -- ++(0, -1) -- ++(7.3, 0) -- (plus-circle.south);
    
        \end{tikzpicture}
        \caption{Pictorial representation of the NN-based bi-fidelity model used as a surrogate for the Godunov flux.}
        \label{fig:nnbifid}
    \end{figure}

   Bi-fidelity learning operates by leveraging the correlation between the data obtained from a high-fidelity (HF) solver, which is usually more accurate but computationally expensive, and a low-fidelity (LF) solver, which is less accurate but relatively inexpensive computationally. Therefore, to capture the nonlinear correlation between the mentioned fidelities, the following approach was put forth by Meng and Karniadakis \cite{meng2020composite}
    \begin{equation}
        y_{\mathrm{H}} = \mathcal{K}(x,y_{\mathrm{L}})
        \label{eqn:arbf1}
    \end{equation}
    where $\mathcal{K}(\cdot)$ represents a potentially nonlinear function, to be determined, that maps between the LF and HF data. Furthermore, $y_\mathrm{L}$, $y_\mathrm{H}$, and $x$ denote the LF, HF, and input data, respectively. 
    
    To approximate the Godunov flux in this case, we utilize a model developed and trained using bi-fidelity learning techniques, which incorporate strategies such as residual learning that model the discrepancy between the HF and LF outputs. The approximate Riemann solver flux, denoted $\Hcal_L^x : \Rbb^m\times\Rbb^m\rightarrow\Rbb^m$, is the LF data with
    \begin{equation}
    	\Hcal_L^x(U^+,U^-) = \Hcal_{roe}(U^+, U^-, e_1), \quad \text{or} \quad
	\Hcal_L^x(U^+,U^-) = \Hcal_{hll}(U^+, U^-, e_1).
    \end{equation}
    The approximation of the one-dimensional Godunov flux using a NN-based bi-fidelity model is
    \begin{equation}
        \mathcal{\hat{H}}_{bf}^x(U^+,U^-; \theta) = \mathcal{H}_L^x(U^+,U^-) + \mathcal{F}_{NN}(U^+,U^-,\mathcal{H}_L^x(U^+,U^-);\theta),
        \label{eqn:nnsm3}
    \end{equation}
    where \( \mathcal{F}_{NN} \) represents the neural network with parameters \( \theta \).
    In this setting, the NN approximates the difference between HF (Godunov flux) and LF (approximate flux) solvers, which can be expressed as the residual
    \begin{equation}
        \mathcal{R}\left(U^+,U^-,\Theta\right) = \Hcal_{gdnv}(U^+,U^-,e_1) - \Theta,
        \label{eqn:nnsm4}
    \end{equation}
    where $\Theta$ is the LF data.
    For training the NN-based bi-fidelity model, we use the following dataset
    \begin{equation}
        \mathcal{D}_{bf} = \left\{ \left( U^+_i, U^-_i, \mathcal{H}_L^x( U^+_i, U^-_i) \right), f_i^{H} = \Hcal_{gdnv}(U^+,U^-,e_1) \right\}_{i=1}^{N_H}.
    \label{eqn:data2}
    \end{equation}    
    The neural network used for the bi-fidelity model is also an FCNN. Furthermore, the pictorial illustration of the bi-fidelity NN (BFNN) model is provided in Figure~\ref{fig:nnbifid}.
    
 

\subsection{Optimization and hyperparameter details}\label{ssec:ohd}

The optimization framework for the NN-based Godunov surrogates is formulated within a supervised learning paradigm. Given the training dataset, the objective is to find the optimal set of parameters $\theta^\star$ that minimizes the discrepancy between the fluxes predicted by the NN-based surrogate and those obtained from the HF solver. Mathematically, this can be expressed as the following optimization problem

\begin{equation}
    \theta^\star = \underset{\theta}{\arg \min } \hspace{0.5em} \mathcal{L}(\mathcal{H}_{gdnv}(U_i^+, U_i^-,e_1), \hat{\mathcal{H}^x}(U_i^+, U_i^-; \theta)),
    \label{eqn:opt}
\end{equation}

\noindent where the surrogate model $\hat{\mathcal{H}}^x$ is either the VNN ($\hat\Hcal_{vn}^x$) or the bi-fidelity model ($\hat\Hcal_{bf}^x$). The goal is to minimize the loss function $\mathcal{L}(\cdot, \cdot)$, which quantifies the above-mentioned discrepancy as

\begin{equation}
    \mathcal{L}(f^H, \hat{f}) = \dfrac{\sum_{i=1}^{N_H} |f_i^H - \hat{f}_i|}{\sum_{i=1}^{N_H} |f_i^H|},
    \label{eqn:loss}
\end{equation}

\noindent where $f^H$ and $\hat{f}$ represent the fluxes from the HF solver and the surrogate model, respectively. The optimization problem (\ref{eqn:opt}) is solved using gradient descent-based algorithms. Specifically, we use Adaptive Moment Estimation (Adam) \cite{kingma2014adam}. In addition, mini-batch training is employed for all surrogate models in this study by subdividing the training set into smaller subsets called mini-batches and using gradient descent for each mini-batch. A suitable schedule is used for the learning rate.
\begin{remark}
 We also evaluated the performance of the proposed models when trained using a relative loss based on $\ell_2$-norm. Across all experiments, it was consistently observed that training with a relative loss based on the $\ell_1$-norm yielded slightly better predictive accuracy compared to the relative loss based on the $\ell_2$-norm.

\end{remark}

\section{Numerical experiments}
\label{sec:num_exp}
In this section,  additional details on the training of the SMs introduced in Sections~\ref{ssec:models:vanilla}-\ref{ssec:models:bifidelity} are provided. Also, the performance of the proposed SMs is studied through a series of a priori and a posteriori tests on one-dimensional and two-dimensional Burgers' and shallow water equations. In the current context, a priori tests directly measure the performance of neural network-based fluxes as surrogates for the Godunov flux. On the other hand, a posteriori analysis refers to evaluating the performance of neural network-based surrogates in downstream tasks for which they were not specifically trained. In this study, the downstream task is an unsteady finite volume method (FVM) simulation that uses the neural network-based surrogates as the numerical flux function. The objective of this analysis is to examine the ability of trained surrogates to generalize beyond the training data and to get insights into their robustness and accuracy as numerical flux functions in unsteady FVM simulations. To save space in the following discussions, the LF solver in the BFNN-based surrogate is built as a correction to the Roe flux; however, similar conclusions are expected to hold when correcting the HLL solver.

\subsection{One-dimensional Burgers' equation}\label{ssec:burg1d}
We begin by considering the nonlinear advection of a scalar quantity $u: \Omega_x\mapsto \mathbb{R}$ through a one-dimensional domain $\Omega_x \subset \mathbb{R}$, which is governed by the following flux functions and source term
\begin{equation} \label{eqn:burg}
    U = u, \quad F_I(U) = \dfrac{u^2}{2}, \quad F_{V}(U,\nabla U) = \nu \nabla u, \quad S(U,\nabla U) = 0,
\end{equation}
\noindent where $\nu: \Omega_x \mapsto \mathbb{R}_{\geq 0}$ is the diffusion coefficient field and $m=1$. We first focus on the construction of NN-based surrogates, i.e., $\hat{\mathcal{H}}^x(U^+,U^-)$, for the inviscid numerical flux. For that, we generate a training dataset by sampling $N_H$ input states, i.e., interior and exterior states, from a uniform distribution $\mathcal{U}(-3.0,3.0)$. For each sampled state, we compute the corresponding HF and LF output data. Once the dataset is created, we train the VNN and BFNN-based SMs using the optimization framework discussed in Section \ref{ssec:ohd}. For the VNN-based SM, we utilize an FCNN comprising three hidden layers, each containing $32$ neurons. Similarly, for the BFNN-based SM, we again use an FCNN with three hidden layers; however, the number of neurons per layer is restricted to $10$. Furthermore, both models are trained for $1300$ epochs. The training process begins with an initial learning rate of $0.01$. A dataset containing $20000$ samples is utilized for training, with a batch size set to $1000$. 

\subsubsection{A priori tests: SMs for inviscid Burgers' flux}
\begin{figure}
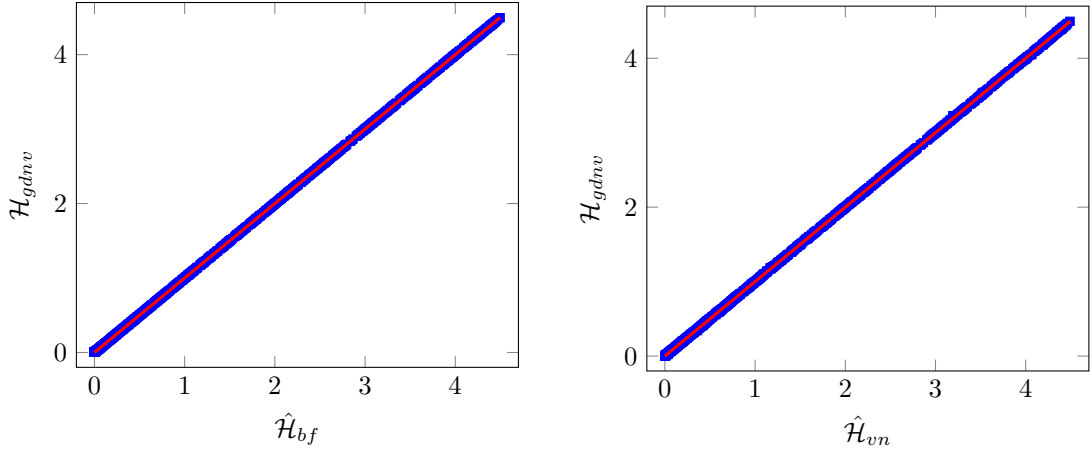

    \centering
    \begin{subfigure}[b]{0.45\textwidth}
        \centering
        \input{_py/burgersdiffapriori.tikz}
        \captionsetup{labelformat=empty}
        \label{fig:aprioriburg11}
    \end{subfigure}
    \begin{subfigure}[b]{0.45\textwidth}
        \centering
        \input{_py/burgersrvnapriori.tikz}
        \captionsetup{labelformat=empty}
        \label{fig:aprioriburg12}
    \end{subfigure}

	\caption{Scatter plots comparing SM-predicted fluxes with the Godunov flux (\ref{scatt:aprbfburg1}) for 1D inviscid Burgers' equation. The best possible match is shown with (\ref{line:aprbfburg1}). Left subfigure: BFNN surrogate vs. Godunov flux and right subfigure: VNN surrogate vs. Godunov flux.}
    \label{fig:aprioriburg1}
\end{figure}

First, we assess the accuracy of the trained SMs by determining how well they predict the Godunov flux on a test set of 2000 samples for the inviscid Burgers' flux in (\ref{eqn:burg}) with $\nu = 0$. The relative $\ell_1$ error over this test set for the BFNN and VNN models are $2.18 \times 10^{-7}$ and $1.4 \times 10^{-3}$, respectively. Also, the error for Roe flux computed on the same test set is $2.1 \times 10^{-1}$. Notably, the test error of the BFNN is four orders of magnitude lower than that of the VNN, despite the latter utilizing a larger network architecture. We also provide the scatter plots that show a comparison between the SM-predicted fluxes and Godunov flux in Figure \ref{fig:aprioriburg1}.

\begin{figure}[h!]
    \centering
    \begin{tikzpicture}
\begin{axis}[
    xlabel={Absolute error},
    ylabel={Frequency},
    xtick={-8.0,-6.0,-4.0,-2.0,0.0},
    xticklabels={$10^{-8}$,$10^{-6}$,$10^{-4}$,$10^{-2}$,$10^{0}$},
    ybar,
    bar width=0.25,
    bar shift=0pt,
]

\addplot[bar width = 0.3,draw=none, fill=red, opacity=1.0] coordinates {
    (-8.0000, 10020)
    (-7.7837, 0)
    (-7.5675, 0)
    (-7.3512, 0)
    (-7.1350, 0)
    (-6.9187, 0)
    (-6.7025, 0)
    (-6.4862, 0)
    (-6.2700, 0)
    (-6.0537, 0)
    (-5.8374, 0)
    (-5.6212, 0)
    (-5.4049, 0)
    (-5.1887, 0)
    (-4.9724, 0)
    (-4.7562, 0)
    (-4.5399, 2)
    (-4.3237, 0)
    (-4.1074, 1)
    (-3.8912, 1)
    (-3.6749, 3)
    (-3.4586, 1)
    (-3.2424, 1)
    (-3.0261, 3)
    (-2.8099, 0)
    (-2.5936, 2)
    (-2.3774, 0)
    (-2.1611, 0)
    (-1.9449, 0)
    (-1.7286, 0)
    (-1.5123, 0)
    (-1.2961, 0)
    (-1.0798, 0)
    (-0.8636, 0)
    (-0.6473, 0)
    (-0.4311, 0)
    (-0.2148, 0)
    (0.0014, 0)
    (0.2177, 0)
    (0.4340, 0)
    (0.6502, 0)
};\label{hist:burg3}
\addplot[draw=none, fill=green, opacity=0.6] coordinates {
    (-8.0000, 9992)
    (-7.7837, 0)
    (-7.5675, 0)
    (-7.3512, 0)
    (-7.1350, 0)
    (-6.9187, 0)
    (-6.7025, 0)
    (-6.4862, 0)
    (-6.2700, 0)
    (-6.0537, 0)
    (-5.8374, 0)
    (-5.6212, 0)
    (-5.4049, 0)
    (-5.1887, 0)
    (-4.9724, 0)
    (-4.7562, 0)
    (-4.5399, 0)
    (-4.3237, 1)
    (-4.1074, 1)
    (-3.8912, 0)
    (-3.6749, 0)
    (-3.4586, 0)
    (-3.2424, 1)
    (-3.0261, 1)
    (-2.8099, 0)
    (-2.5936, 1)
    (-2.3774, 0)
    (-2.1611, 2)
    (-1.9449, 0)
    (-1.7286, 1)
    (-1.5123, 0)
    (-1.2961, 0)
    (-1.0798, 0)
    (-0.8636, 0)
    (-0.6473, 0)
    (-0.4311, 0)
    (-0.2148, 0)
    (0.0014, 0)
    (0.2177, 0)
    (0.4340, 0)
    (0.6502, 0)
};\label{hist:burg2}
\addplot[bar width = 0.2,draw=none, fill=blue, opacity=0.9] coordinates {
    (-8.0000, 1)
    (-7.7837, 1)
    (-7.5675, 0)
    (-7.3512, 0)
    (-7.1350, 0)
    (-6.9187, 0)
    (-6.7025, 1)
    (-6.4862, 0)
    (-6.2700, 2)
    (-6.0537, 3)
    (-5.8374, 1)
    (-5.6212, 4)
    (-5.4049, 7)
    (-5.1887, 8)
    (-4.9724, 10)
    (-4.7562, 10)
    (-4.5399, 9)
    (-4.3237, 23)
    (-4.1074, 25)
    (-3.8912, 30)
    (-3.6749, 39)
    (-3.4586, 50)
    (-3.2424, 58)
    (-3.0261, 86)
    (-2.8099, 86)
    (-2.5936, 115)
    (-2.3774, 206)
    (-2.1611, 194)
    (-1.9449, 278)
    (-1.7286, 337)
    (-1.5123, 417)
    (-1.2961, 539)
    (-1.0798, 695)
    (-0.8636, 830)
    (-0.6473, 912)
    (-0.4311, 1085)
    (-0.2148, 1175)
    (0.0014, 1240)
    (0.2177, 1022)
    (0.4340, 501)
    (0.6502, 0)
};\label{hist:burg1}

\end{axis}
\end{tikzpicture}
    \caption{Histogram of absolute error between the Godunov flux outputs and the predictions from the BFNN {\large(}\ref{hist:burg2}{\large)}, VNN {\large(}\ref{hist:burg3}{\large)}, and Roe {\large(}\ref{hist:burg1}{\large)} flux for the a priori rarefaction test.}
    \label{fig:aprioriburg2}
\end{figure}
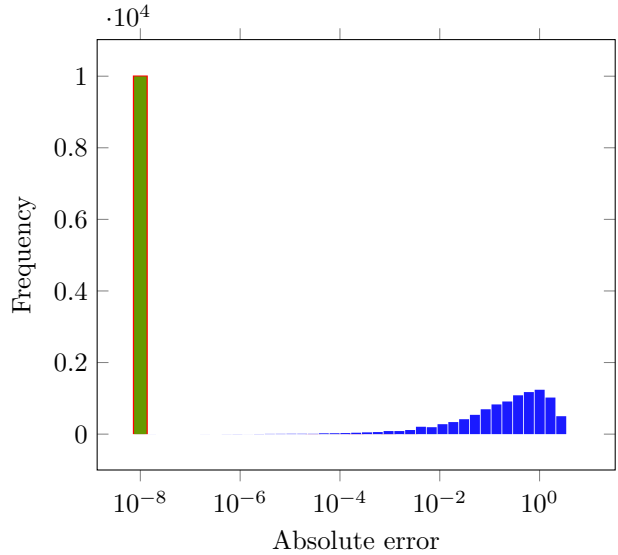

Moreover, as discussed in Section \ref{sec:intro}, approximate Riemann solvers encounter difficulties in handling some special Riemann problem cases. Therefore, we create a challenging a priori test case and compare the performance of trained surrogate models and Roe flux to the Godunov flux. We create the test set by generating $10000$ input samples such that $U^+$ and $U^-$ are sampled from uniform distributions $\mathcal{U}(0.0,3.0)$  and $\mathcal{U}(-3.0,0.0)$, respectively, which ensures the solution will be a rarefaction wave whose flux as $x/t=0$ is zero. The comparison between the predictions generated by the trained SMs and the outputs of the Roe flux function, in terms of absolute error computed against the corresponding outputs of the Godunov flux function, is presented in Figure \ref{fig:aprioriburg2}. Figure \ref{fig:aprioriburg2} indicates that the predictions from both surrogate models exhibit consistently lower errors compared to the Roe flux.

\subsubsection{A posteriori test: Inviscid Burgers' equation in one dimension}
\begin{figure}
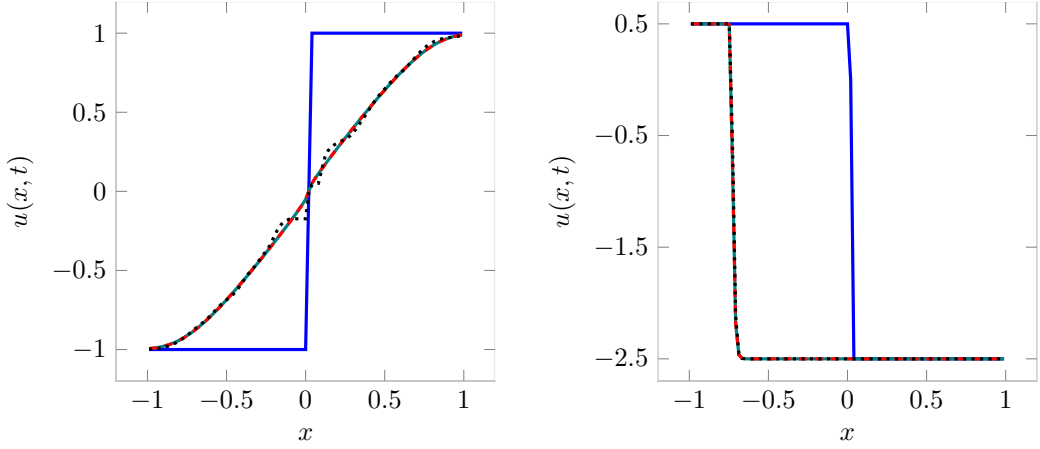

    \centering
    \begin{subfigure}[b]{0.4\textwidth}
        \centering
        \input{_py/burgersfvpost1.tikz}
        \captionsetup{labelformat=empty}
        \label{fig:apostburg11}
    \end{subfigure}
    \hspace{0.02\textwidth} 
    \begin{subfigure}[b]{0.4\textwidth}
        \centering
        \input{_py/burgersfvpost2.tikz}
        \captionsetup{labelformat=empty}
        \label{fig:apostburg12}
    \end{subfigure}

    \caption{FVM-based numerical solutions for the 1D inviscid Burgers equation at time $t = 0.75\,s$. Solutions are obtained using numerical fluxes from BFNN (\ref{line:burg1fv13}), VNN (\ref{line:burg1fv14}), and the Godunov (\ref{line:burg1fv12}) solvers for two distinct initial conditions (\ref{line:burg1fv11}). Left subfigure: Case-I and right subfigure: Case-II.}
    \label{fig:apostburg1}
\end{figure}

In this section, we evaluate the trained SMs as numerical flux functions for an FVM simulation. The assessment focuses on solving the one-dimensional inviscid Burgers equation, (\ref{eqn:burg}) with $\nu = 0$, with an FVM-based solver on a spatial domain $\Omega_{x}:= (-1, 1)$ over a time interval $\mathcal{T} = (0,0.75)$ for the two initial conditions

\begin{equation*}
  \textbf{Case I: } u(x, 0) = 
\begin{cases} 
-1 & x \in [-1.0, 0.0), \\
 1 & x \in [0.0, 1.0].
\end{cases}
\quad 
\quad \textbf{Case II: } u(x, 0) = 
\begin{cases} 
0.5 & x \in [-1,0, 0.0), \\
\color{red}{-2.5} & x \in [0.0, 1].
\end{cases}  
\end{equation*}
For the given initial conditions, the solutions at the end of the time interval of interest $(t = 0.75\,s)$, obtained using the trained SMs as numerical flux functions, are presented in Figure \ref{fig:apostburg1}. For comparison, the corresponding Godunov flux-based solutions are also plotted in the same figure. As shown in Figure \ref{fig:apostburg1}, the VNN flux-based solution exhibits oscillations in the solution at the time of interest. In contrast, the BFNN flux-based solution closely aligns with the Godunov flux-based solution. This comparison highlights the superior generalization capability of the BFNN-based SM over the VNN-based SM for downstream tasks.  

\subsubsection{A posteriori test: The viscous Burgers' equation in one dimension}
\begin{figure}
    \centering

    \begin{center}
        \begin{subfigure}[b]{0.45\textwidth}
            \centering
            \input{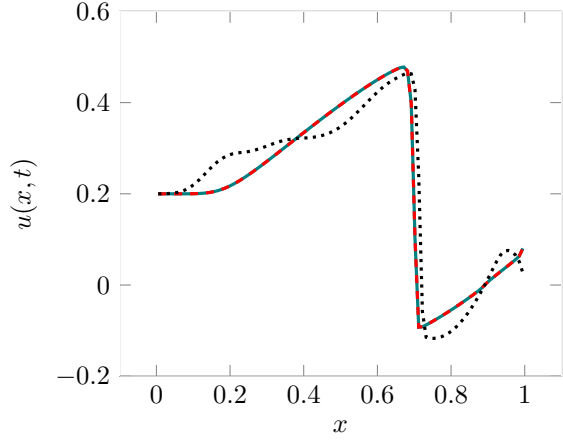}
        \end{subfigure}
    \end{center}
    \caption{FVM-based numerical solutions for the one-dimensional viscous Burgers equation at time $t = 1.0 \,s$. Solutions are obtained using numerical fluxes from BFNN (\ref{line:burg1fv13}), VNN (\ref{line:burg1fv14}), and the Godunov (\ref{line:burg1fv12}) solvers.}
    \label{fig:apostvburg2}
\end{figure}

\begin{figure}
    \centering

    \begin{center}
        \begin{subfigure}[b]{0.45\textwidth}
            \centering
            \input{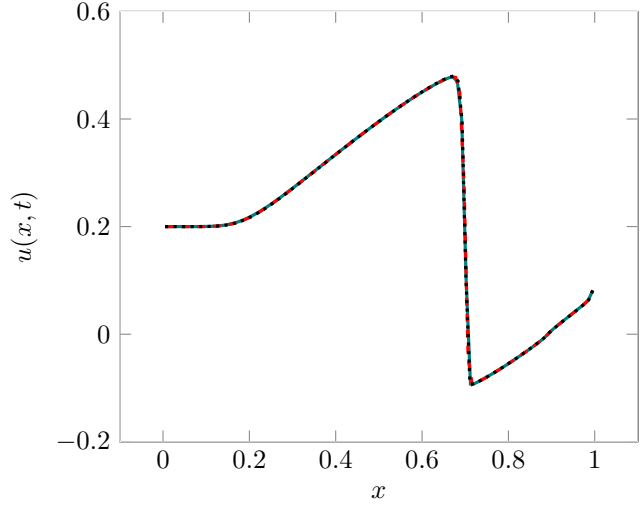}
        \end{subfigure}
    \end{center}
    \caption{FVM-based numerical solutions for the one-dimensional viscous Burgers equation at time $t = 1.0 \,s$ for grid spacing $\Delta x = 1.01\times 10^{-2}$ (\ref{line:burg4fv11}), $\Delta x = 9.60 \times 10^{-3}$ (\ref{line:burg4fv12}), and $\Delta x = 9.20 \times 10^{-3}$ (\ref{line:burg4fv13}) using BFNN as inviscid numerical flux function.}
    \label{fig:apostvburg3}
\end{figure}

Fluxes obtained from Riemann solvers also play a key role in numerical simulations of viscous flows, particularly those dominated by convection.
Therefore, we subject the trained SMs to further evaluation as inviscid numerical flux functions in an FVM setting by solving viscous Burgers' equation (\ref{eqn:burg}) with constant diffusion coefficient $\nu = 10^{-4}$ in the spatial domain $\Omega_{x}:= (0, 1)$ and temporal domain $\mathcal{T} = (0,1.0)$. Note that the viscous flux is computed using a second-order FVM scheme. Furthermore, we choose the smooth initial condition
\begin{equation}
    u(x,0) = 0.2 + \frac{\sin(2\pi x)}{\pi}.
\end{equation}
Figure \ref{fig:apostvburg2} presents the FVM solutions obtained using different SM-based flux functions at time $t = 1.0\,s$, alongside the reference solution computed from the Godunov solver. The figure illustrates that the simulation conducted with the VNN-based flux function produces a solution that deviates from the Godunov reference. In contrast, the solution obtained using the BFNN-based flux function closely aligns with the Godunov solution. This observation further underscores the relatively better generalization capabilities of the BFNN flux framework compared to the VNN-based approach. To substantiate the proposed models' compatibility with multiple grid resolutions, a characteristic shared with the Godunov flux, FVM simulations were executed on different computational grids using BFNN as the inviscid numerical flux function. The resulting solution profiles at $t = 1.0\,s$, depicted in Figure \ref{fig:apostvburg3}, corroborate our previous claim. Similar conclusions hold for VNN, as its inputs also comprise quantities that are mesh-agnostic.

\subsection{One-dimensional shallow water equations}\label{ssec:swe1d}
Next, we focus on the shallow water equations (SWE), which describe the behavior of an incompressible, inviscid fluid under the influence of gravity. These equations assume no vertical acceleration and neglect heat conduction. The analysis is restricted to a one-dimensional domain, $\Omega_x \subset \mathbb{R}$. In this context, the SWE in one dimension can be expressed as
\begin{equation}
    U = \begin{bmatrix} h \\ hu \end{bmatrix}, \quad F_I(U) = \begin{bmatrix} hu \\
    hu^2 + \frac{1}{2}gh^2 \end{bmatrix}, \quad F_{V}(U,\nabla U) = 0, \quad S(U,\nabla U) = 0,
    \label{eqn:swe1d}
\end{equation}
\noindent where $h(x,t) : \Omega_x \to \mathbb{R}_{>0}$ is the height of the free surface, $u(x,t): \Omega_x \to \mathbb{R}$ is the fluid velocity, and $g \in \mathbb{R}_{>0}$ is the gravitational acceleration. In this study, we $g$ is set to be equal to one. The procedure for constructing the surrogate is similar to the approach used for the inviscid Burgers' equation outlined in Section \ref{ssec:burg1d}. In this case, however, the input states, i.e., $U^+$ and $U^-$, are generated by sampling primitive variables, height $(h)$ and velocity $(u)$, from the uniform distributions $\mathcal{U}(0,3.5)$ and $\mathcal{U}(-2.5,2.5)$, respectively. FCNNs with four hidden layers and $40$ neurons per layer are used for both VNN and BFNN surrogates. Both models are trained over a total of $1500$ epochs. The training process is conducted with an initial learning rate of $0.01$. A dataset comprising $45000$ samples is used for training, with a batch size set to $3000$. 

\subsubsection{A priori tests: SMs for inviscid SWE flux}
\begin{figure}
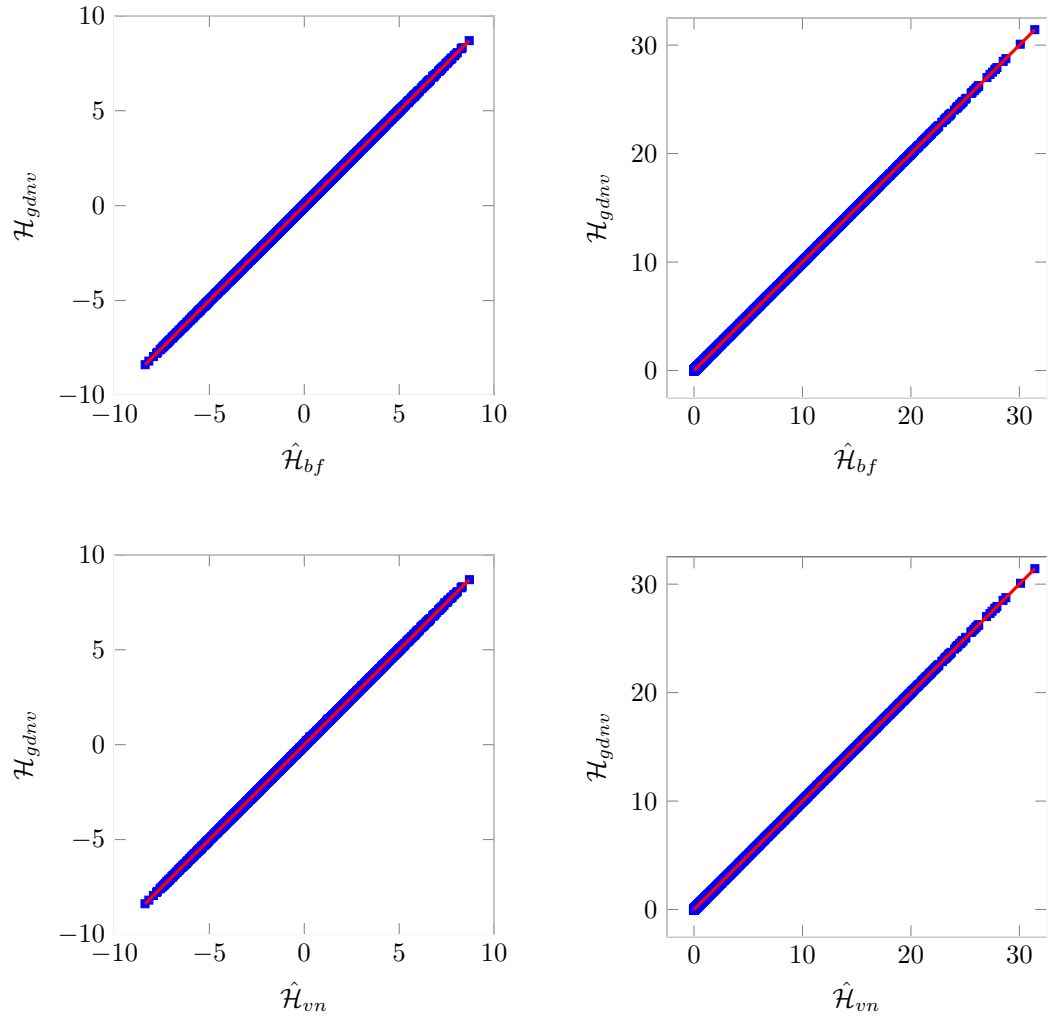

    \centering
    \begin{minipage}{\textwidth}
        \centering
        \begin{subfigure}[t]{0.4\textwidth}
            \centering
            \input{_py/swebfapriori_1stcomp.tikz}
            \captionsetup{labelformat=empty}
            \label{fig:aprioriswe11}
        \end{subfigure}
        \hspace{0.04\textwidth}
        \begin{subfigure}[t]{0.4\textwidth}
            \centering
            \input{_py/swebfapriori_2ndcomp.tikz}
            \captionsetup{labelformat=empty}
            \label{fig:aprioriswe13}
        \end{subfigure}
    \end{minipage}

    \vspace{0.25cm}

    \begin{minipage}{\textwidth}
        \centering
        \begin{subfigure}[t]{0.4\textwidth}
            \centering
            \input{_py/swevnapriori_1stcomp.tikz}
            \captionsetup{labelformat=empty}
            \label{fig:aprioriswe12}
        \end{subfigure}
        \hspace{0.04\textwidth} 
        \begin{subfigure}[t]{0.4\textwidth}
            \centering
            \input{_py/swevnapriori_2ndcomp.tikz}
            \captionsetup{labelformat=empty}
            \label{fig:aprioriswe14}
        \end{subfigure}
    \end{minipage}

	\caption{Scatter plots comparing the BFNN (\textit{top}) and VNN (\textit{bottom}) fluxes with the Godunov flux for one-dimensional SWE. The \textit{left} column shows the first component of the flux ($hu$) and the \textit{right} column shows the second component ($hu^2+gh^2/2$). The best possible match is shown with (\ref{line:aprbfburg1}).}
    \label{fig:aprioriswe1}
\end{figure}

To evaluate the performance of the trained SMs for SWE flux, the relative $\ell_1$ errors were calculated on a test set containing $5000$ samples. The errors for the BFNN and VNN models are $6.0\times 10^{-4}$ and $1.3\times 10^{-3}$, respectively. Furthermore, the error for Roe flux outputs, used as inputs in the BFNN-based surrogate, with respect to the Godunov flux outputs was found to be $7.2 \times 10^{-2}$. A visual representation of the trained models' capabilities for flux prediction is provided through scatter plots in Figure \ref{fig:aprioriswe1}, which depict a comparison between the predicted flux components from the trained models and those from the Godunov solver.

\begin{figure}
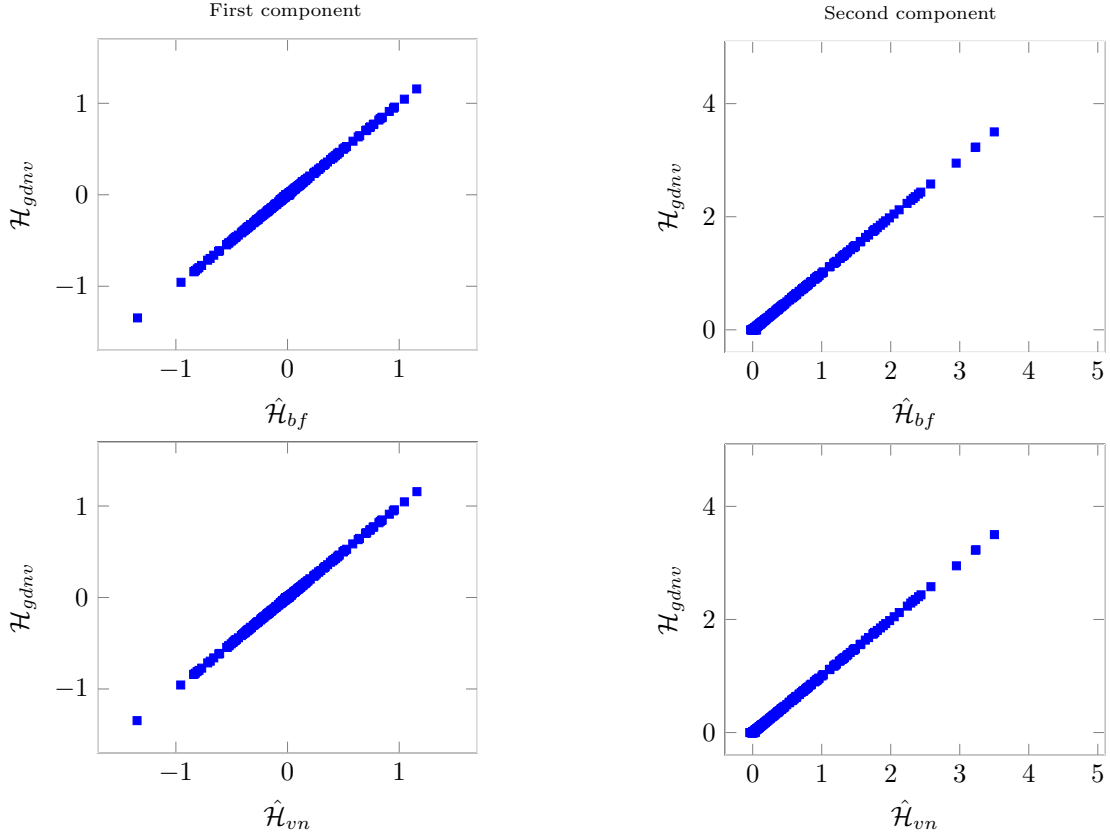

    \centering
    \begin{minipage}{1.0\textwidth}
        \centering
        \begin{subfigure}[b]{0.4\textwidth}
            \centering
            \captionsetup{labelformat=empty, position=top, justification=centering}
            \caption{\hspace{1cm} First component}
            \input{_py/swetrs1bf_1stcomp.tikz}
        \end{subfigure}
        \hspace{0.1\textwidth}
        \begin{subfigure}[b]{0.4\textwidth}
            \centering
            \captionsetup{labelformat=empty, position=top, justification=centering}
            \caption{\hspace{2em} Second component}
            \input{_py/swetrs1bf_2ndcomp.tikz}
        \end{subfigure}
    \end{minipage}
    \vspace{0.4cm}
    \begin{minipage}{1.0\textwidth}
        \centering
        \begin{subfigure}[b]{0.4\textwidth}
            \centering
            \input{_py/swetrs1vn_1stcomp.tikz}
        \end{subfigure}
        \hspace{0.1\textwidth}
        \begin{subfigure}[b]{0.4\textwidth}
            \centering
            \input{_py/swetrs1vn_2ndcomp.tikz}
        \end{subfigure}
    \end{minipage}
    \caption{Scatter plots comparing the BFNN (\textit{first row}) and VNN (\textit{second row}) fluxes with Godunov flux for the one-dimensional SWE (scenario-I). The \textit{left} column shows the first component of the flux ($hu$) and the \textit{right} column shows the second component ($hu^2+gh^2/2$).}
    \label{fig:combined_aprioriswe_1}
\end{figure}

\begin{figure}
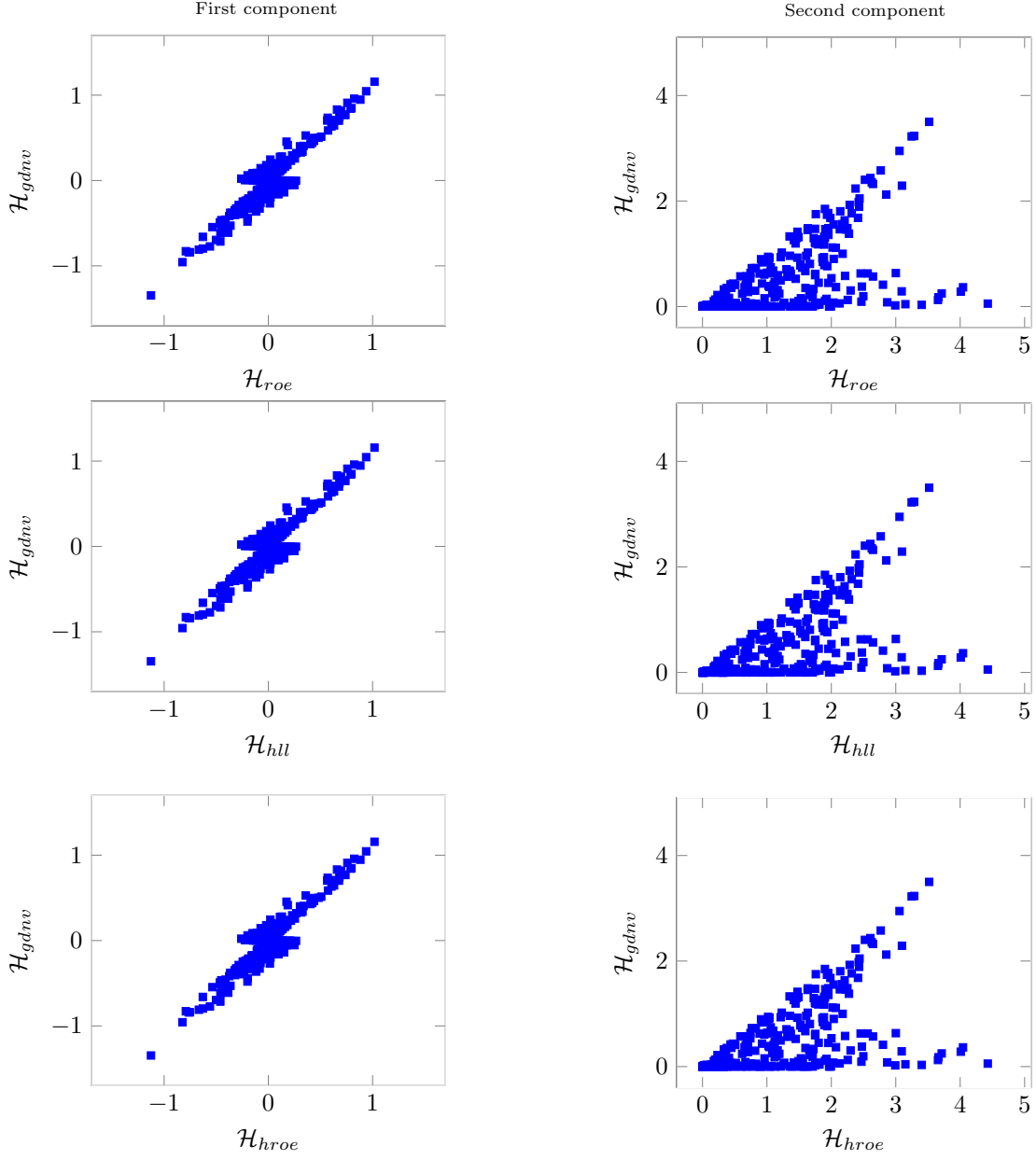

    \centering
    \begin{minipage}{1.0\textwidth}
        \centering
        \begin{subfigure}[b]{0.4\textwidth}
            \centering
            \captionsetup{labelformat=empty, position=top, justification=centering}
            \caption{\hspace{1cm} First component}
            \input{_py/swetrs1roe_1stcomp.tikz}
        \end{subfigure}
        \hspace{0.1\textwidth}
        \begin{subfigure}[b]{0.4\textwidth}
            \centering
            \captionsetup{labelformat=empty, position=top, justification=centering}
            \caption{\hspace{1cm} Second component}
            \input{_py/swetrs1roe_2ndcomp.tikz}
        \end{subfigure}
    \end{minipage}
    \vspace{0.4cm}
    \begin{minipage}{1.0\textwidth}
        \centering
        \begin{subfigure}[b]{0.4\textwidth}
            \centering
            \input{_py/swetrs1hll_1stcomp.tikz}
        \end{subfigure}
        \hspace{0.1\textwidth}
        \begin{subfigure}[b]{0.4\textwidth}
            \centering
            \input{_py/swetrs1hll_2ndcomp.tikz}
        \end{subfigure}
    \end{minipage}
    \vspace{0.4cm}
    \begin{minipage}{1.0\textwidth}
        \centering
        \begin{subfigure}[b]{0.4\textwidth}
            \centering
            \input{_py/swetrs1roef_1stcomp.tikz}
        \end{subfigure}
        \hspace{0.1\textwidth}
        \begin{subfigure}[b]{0.4\textwidth}
            \centering
            \input{_py/swetrs1roef_2ndcomp.tikz}
        \end{subfigure}
    \end{minipage}
    \caption{Scatter plots comparing Roe (\textit{first row}), HLL (\textit{second row}), and Roe with Harten fix (\textit{third row}) fluxes with Godunov flux for the one-dimensional SWE (scenario-I). The \textit{left} column shows the first component of the flux ($hu$) and the \textit{right} column shows the second component ($hu^2+gh^2/2$).}
    \label{fig:combined_aprioriswe_2}
\end{figure}

In addition, a test scenario that presents a significant challenge to approximate Riemann solvers is defined as 
\begin{equation*}
\textbf{Scenario-I:}  \qquad
\begin{cases} 
h^- = h^+,\, 0<h^+<3, \\
-2<u^-<0,\, 0<u^+<2.
\end{cases}
\end{equation*}
\noindent These constraints lead to a Riemann problem whose solution contains a rarefaction wave and a height that diminishes with time. The flux outputs for scenario-I obtained from Roe, HLL, and Roe solver with Harten fix along with the flux predictions from the trained SMs are plotted against the Godunov flux outputs in Figures \ref{fig:combined_aprioriswe_1} and \ref{fig:combined_aprioriswe_2}. These figures show that, for scenario-I, the flux predictions from the NN-based surrogates exhibit negligible deviation from the target Godunov flux. This deviation is notably smaller compared to fluxes obtained from approximate Riemann solvers. 

\subsubsection{A posteriori test: SWE in one dimension}
We now conduct FVM simulations of the one-dimensional SWE (\ref{eqn:swe1d}) using trained SMs as the numerical flux functions. The spatial domain is $\Omega_x := (-1,1)$ and the simulations are run for the following three initial conditions

\begin{align*}
\textbf{Case-I: } &\quad h(x, 0) =  1, \quad x \in [-1.0, 1.0],
&\quad u(x, 0) = 
\begin{cases} 
-0.5 & x \in [-1.0, 0.0), \\
0.5 & x \in [0.0, 1.0].
\end{cases}  \\[1em]
\textbf{Case-II: } &\quad h(x, 0) =
\begin{cases} 
2.0 & x \in [-1.0, 0.0), \\
1.0 & x \in [0.0, 1.0].
\end{cases},
&\quad u(x, 0) = 
\begin{cases} 
-0.5 & x \in [-1.0, 0.0), \\
0.5 & x \in [0.0, 1.0].
\end{cases}\\[1em]  
\textbf{Case-III: } &\quad h(x, 0) =  1, \quad x \in [-1.0, 1.0]
&\quad u(x, 0) = 
\begin{cases} 
-1.0 & x \in [-1.0, 0.0), \\
1.0 & x \in [0.0, 1.0].
\end{cases}  
\end{align*}

\begin{figure}
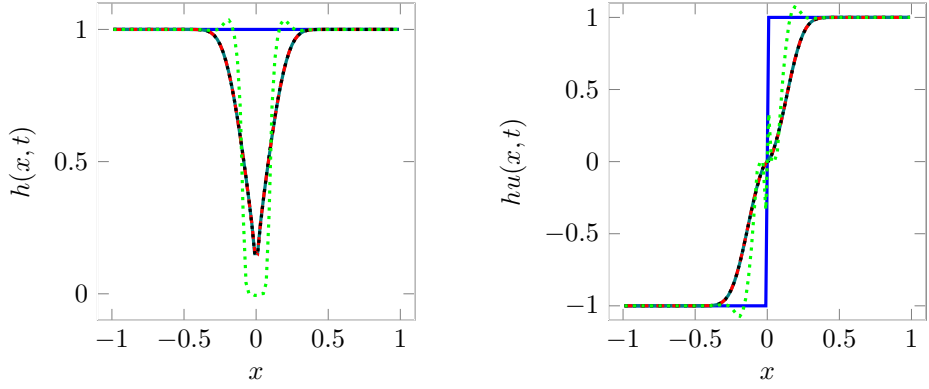

    \centering
    \begin{subfigure}[t]{0.35\textwidth}
        \centering
        \input{_py/fail_swefvcomp1post3.tikz}
        \captionsetup{labelformat=empty}
        \label{fig:apostswe31}
    \end{subfigure}
    \hspace{0.03\textwidth} 
    \begin{subfigure}[t]{0.35\textwidth}
        \centering
        \input{_py/fail_swefvcomp2post3.tikz}
        \captionsetup{labelformat=empty}
        \label{fig:apostswe32}
    \end{subfigure}

    \caption{FVM-based numerical solutions for the one-dimensional SWEs at time $t = 0.1\,s$. Solutions are obtained using numerical fluxes from BFNN (\ref{line:swe3fv13}), VNN (\ref{line:swe3fv14}), HLL (\ref{line:swe3fv15}), and the Godunov solver (\ref{line:swe1fv12}) for initial condition (\ref{line:swe3fv11}) in case-III.}
    \label{fig:apostswe3}
\end{figure}

\noindent For all three cases, the left and right states are connected through rarefactions. Additionally, the initial condition evolves into a near-dry state ($h=0$) over time for case-III. It is well established that the HLL solver is more robust than the Roe solver in handling near-dry states \cite{ketcheson2020riemann}.
The FVM solutions obtained using different SM-based flux functions for case-I and case-II are presented in Figure \ref{fig:apostswe1} and \ref{fig:apostswe2}, along with the reference solution obtained using Godunov flux-based FVM simulation. 
The solutions obtained for case-III using SM-based fluxes, the HLL solver, and the Godunov flux function are presented in Figure \ref{fig:apostswe3}. Simulations were also attempted for this case using the Roe solver and the Roe solver with the Harten fix. However, these simulations failed and terminated prematurely, well before the final time of interest. As a result, the corresponding solutions are not included in Figure \ref{fig:apostswe3}. 
Furthermore, it can be observed in Figure \ref{fig:apostswe3} that even the solutions provided by the HLL solver for both SWE components do not match the Godunov-based ground truth. In contrast, the solutions obtained using SM-based fluxes closely follow the Godunov-based reference.

\begin{figure}
    \centering
    \begin{subfigure}[t]{0.35\textwidth}
        \centering
        \input{_py/swefvcomp1post2.tikz}
        \captionsetup{labelformat=empty}
        \label{fig:apostswe11}
    \end{subfigure}
    \hspace{0.03\textwidth} 
    \begin{subfigure}[t]{0.35\textwidth}
        \centering
        \input{_py/swefvcomp2post2.tikz}
        \captionsetup{labelformat=empty}
        \label{fig:apostswe12}
    \end{subfigure}
    \vspace{-1em}
    \caption{FVM-based numerical solutions for the one-dimensional SWEs at time $t = 0.1\,s$. Solutions are obtained using numerical fluxes from BFNN (\ref{line:swe1fv13}), VNN (\ref{line:swe1fv14}), and the Godunov solver (\ref{line:swe1fv12}) for the initial condition (\ref{line:swe1fv11}) in case-I.}
    \label{fig:apostswe1}
\end{figure}
\begin{figure}
    \centering 
    \begin{subfigure}[t]{0.35\textwidth}
        \centering
        \input{_py/swefvcomp1post1.tikz}
        \captionsetup{labelformat=empty}
        \label{fig:apostswe21}
    \end{subfigure}
    \hspace{0.03\textwidth} 
    \begin{subfigure}[t]{0.35\textwidth}
        \centering
        \input{_py/swefvcomp2post1.tikz}
        \captionsetup{labelformat=empty}
        \label{fig:apostswe22}
    \end{subfigure}
    \caption{FVM-based numerical solutions for the one-dimensional SWEs at time $t = 0.2\,s$. Solutions are obtained using numerical fluxes from BFNN (\ref{line:swe1fv13}), VNN (\ref{line:swe1fv14}), and the Godunov solver (\ref{line:swe1fv12}) for the initial condition (\ref{line:swe1fv11}) in case-II.}
    \label{fig:apostswe2}
\end{figure}

\subsection{Two-dimensional inviscid Burgers' equation} \label{ssec:burg2d}
In this example, we consider nonlinear advection of a scalar quantity 
$u: \Omega_x \mapsto \mathbb{R}$ within a two-dimensional domain 
$\Omega_x \subset \mathbb{R}^2$. The governing equations are defined by the flux functions and source term below
\begin{equation}
    U = u, \quad F_I(U) = \dfrac{u^2}{2}\beta^T, \quad F_{V}(U,\nabla U) = 0, \quad S(U,\nabla U) = 0,
    \label{eqn:burg2d}
\end{equation}
where $\beta = (\beta_1, \beta_2)$ is the advection velocity with $\beta_i:\Omega_x \mapsto \mathbb{R}$. In this problem, we take the advection field to be constant with $\beta_1 = \beta_2 = 1$. In this multi-dimensional setting, we do not develop a surrogate for the multi-dimensional numerical flux; rather, we reuse the one-dimensional flux surrogates developed in the previous section by exploiting the symmetries in the governing equations (Section \ref{sec:prop}). 

\begin{figure}
    \centering
    \begin{minipage}{\textwidth}
        \centering
        \begin{subfigure}[b]{0.4\textwidth}
            \centering
            \begin{tikzpicture}
\begin{axis}[
axis equal image,
xmin=-0.9416904449462891,
xmax=0.9417643547058105,
ymin=-0.8028290867805481,
ymax=0.987214982509613,
width=1.0\textwidth,
xtick={-.9, 0, 0.9},
ytick={-.8,0, 0.9},
width=1.\textwidth]
\addplot []
graphics [xmin=-0.9416904449462891,xmax=0.9417643547058105,ymin=-0.8028290867805481,ymax=0.987214982509613] {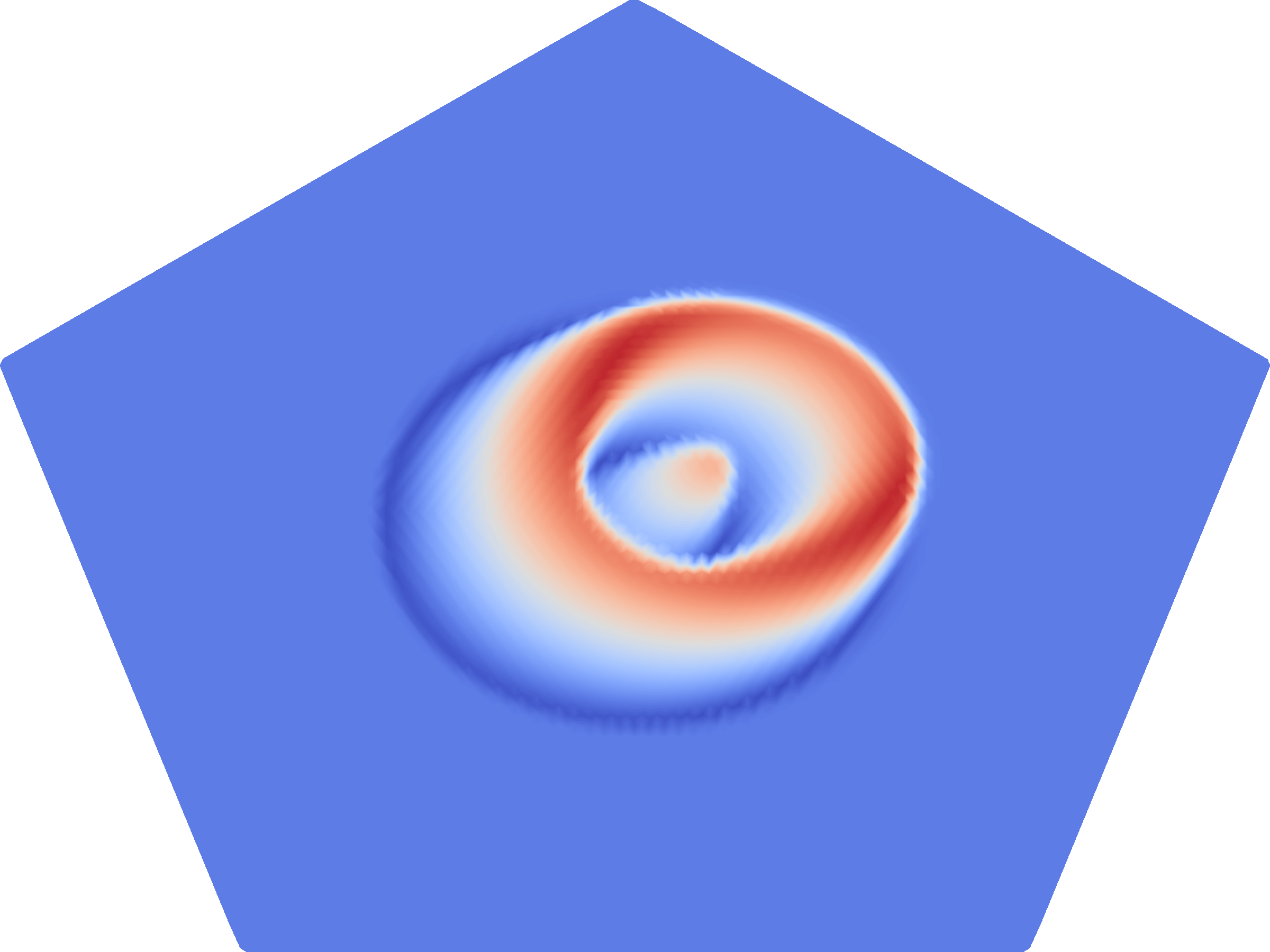};

\end{axis}
\end{tikzpicture}
            \\
            \hspace{2em}\colorbarMatplotlibCoolwarm{0}{0.25}{0.50}{0.75}{1.0}{5cm}{0.3cm} 
            \captionsetup{labelformat=empty}
            \caption{\hspace{0.6cm}(a) Solution at time $t = 0.2\,s$}
            \label{fig:apostbe21}
        \end{subfigure}
        \hspace{0.01\textwidth}
        \begin{subfigure}[b]{0.4\textwidth}
            \centering
            \begin{tikzpicture}
\begin{axis}[
axis equal image,
xmin=-0.9416904449462891,
xmax=0.9417643547058105,
ymin=-0.8028290867805481,
ymax=0.987214982509613,
width=1.0\textwidth,
xtick={-.9, 0, 0.9},
ytick={-.8, 0, 0.9},
width=1.\textwidth]
\addplot []
graphics [xmin=-0.9416904449462891,xmax=0.9417643547058105,ymin=-0.8028290867805481,ymax=0.987214982509613] {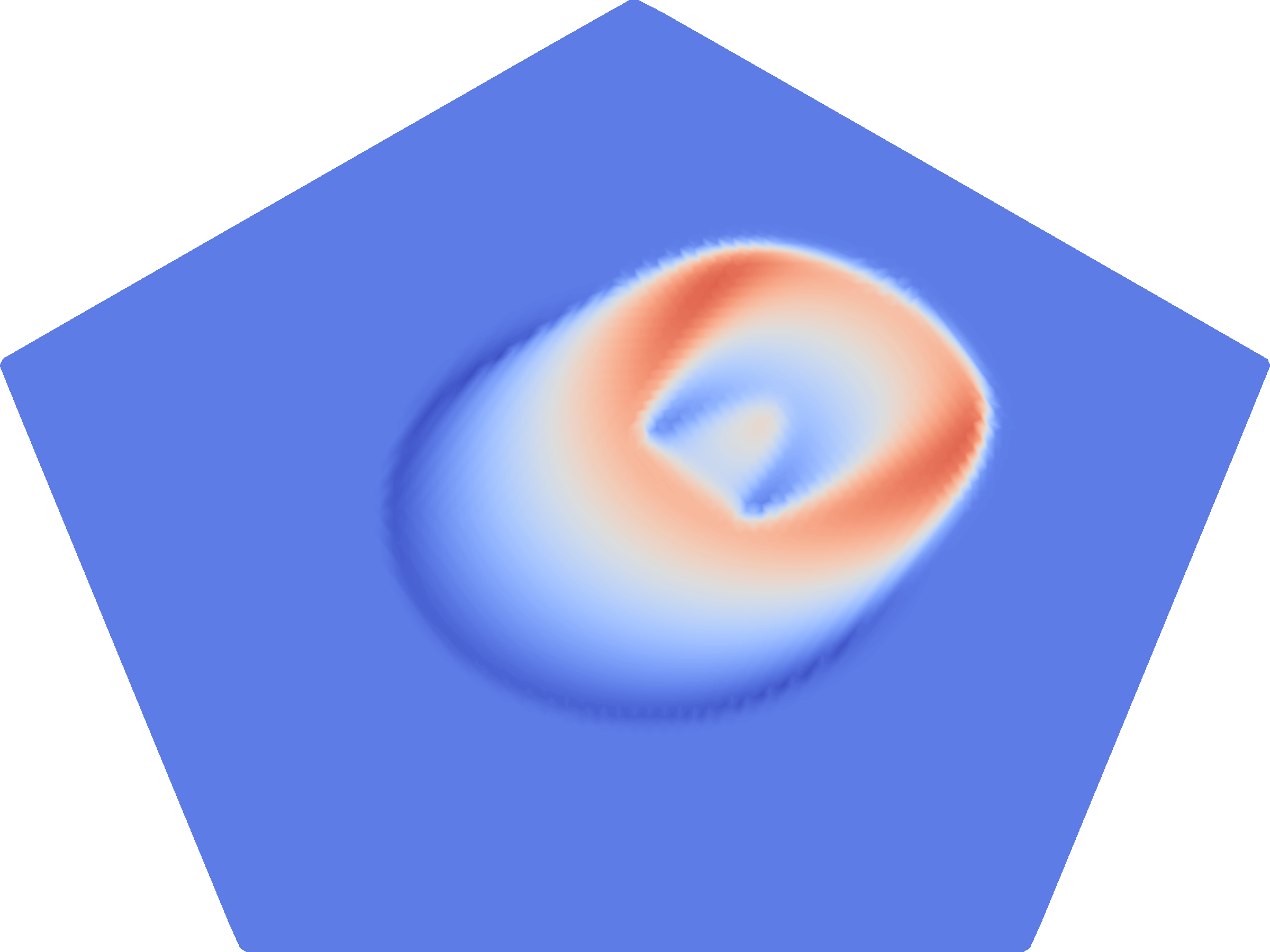};

\end{axis}
\end{tikzpicture}
            \\
            \hspace{2em}\colorbarMatplotlibCoolwarm{0}{0.25}{0.50}{0.75}{1.0}{5cm}{0.3cm} 
            \captionsetup{labelformat=empty}
            \caption{\hspace{0.7cm}(b) Solution at time $t = 0.4\,s$}
            \label{fig:apostbe22}
        \end{subfigure}
    \end{minipage}
    
    \vspace{1cm}
    
    \begin{minipage}{\textwidth}
        \centering
        \begin{subfigure}[b]{0.4\textwidth}
            \centering
            \begin{tikzpicture}
\begin{axis}[
axis equal image,
xmin=-0.9416904449462891,
xmax=0.9417643547058105,
ymin=-0.8028290867805481,
ymax=0.987214982509613,
xtick={-.9,0,0.9},
ytick={-.8,0,0.9},
width=1.0\textwidth]
\addplot []
graphics [xmin=-0.9416904449462891,xmax=0.9417643547058105,ymin=-0.8028290867805481,ymax=0.987214982509613] {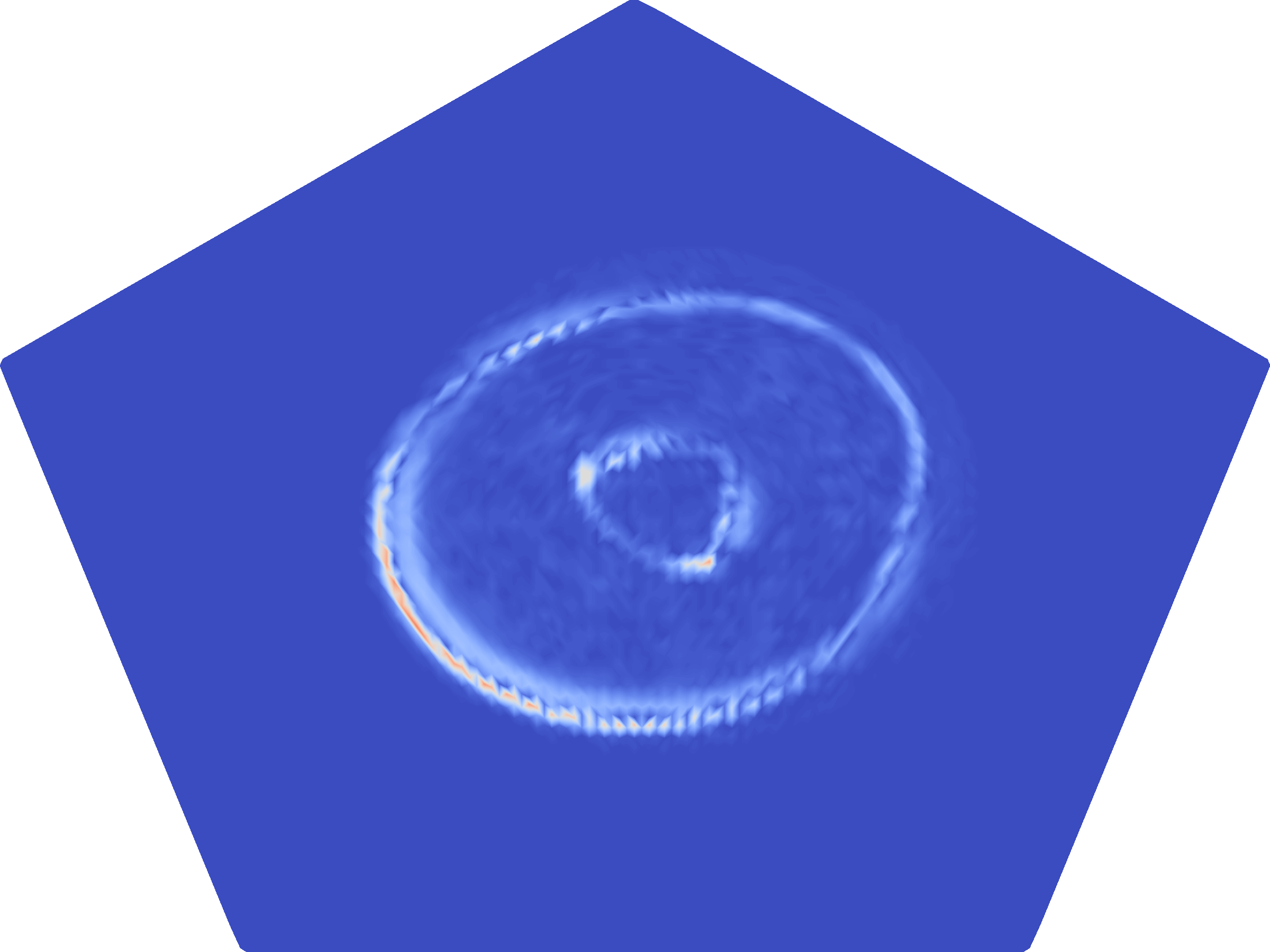};

\end{axis}
\end{tikzpicture}
            \\
            \hspace{2em}\colorbarMatplotlibCoolwarm{0}{1.0e-6}{2.0e-6}{3.0e-6}{4.0e-6}{5cm}{0.3cm} 
            \captionsetup{labelformat=empty}
            \caption{\hspace{0.85cm}(c) Absolute error at time $t = 0.2\,s$}
            \label{fig:apostbe23}
        \end{subfigure}
        \hspace{0.01\textwidth}
        \begin{subfigure}[b]{0.4\textwidth}
            \centering
            \begin{tikzpicture}
\begin{axis}[
axis equal image,
xmin=-0.9416904449462891,
xmax=0.9417643547058105,
ymin=-0.8028290867805481,
ymax=0.987214982509613,
xtick={-.9,0,0.9},
ytick={-.8, 0,0.9},
width=1.0\textwidth]
\addplot []
graphics [xmin=-0.9416904449462891,xmax=0.9417643547058105,ymin=-0.8028290867805481,ymax=0.987214982509613] {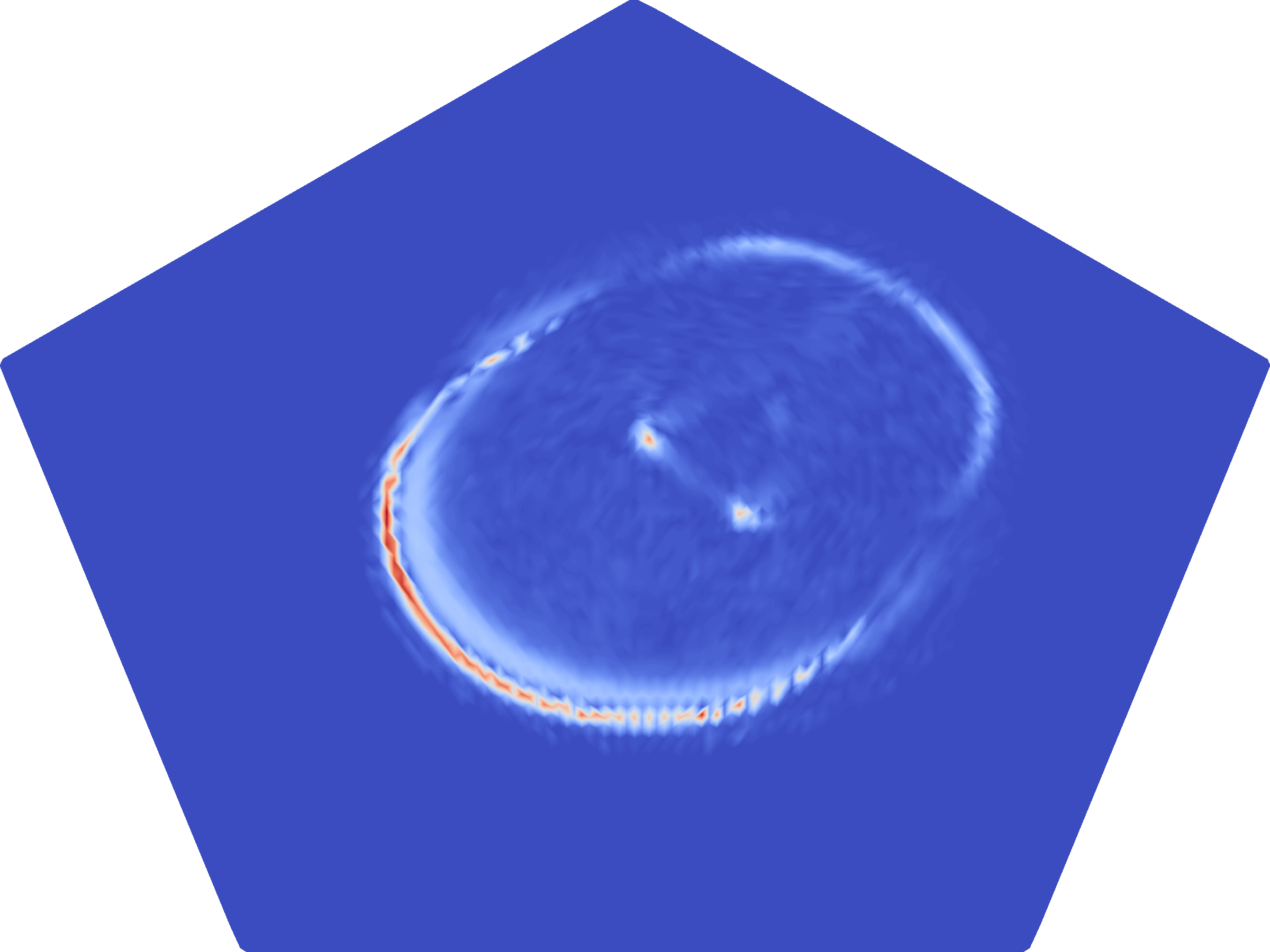};

\end{axis}
\end{tikzpicture}
            \\
            \hspace{2em}\colorbarMatplotlibCoolwarm{0}{1.0e-6}{2.0e-6}{3.0e-6}{4.0e-6}{5cm}{0.3cm} 
            \captionsetup{labelformat=empty}
            \caption{\hspace{0.59cm}(d) Absolute error at time $t = 0.4\,s$}
            \label{fig:apostbe24}
        \end{subfigure}
    \end{minipage}

	\caption{FVM-based numerical solutions obtained using BFNN as numerical flux at different time instances and corresponding absolute error field (relative to FVM simulation with Godunov numerical flux) for two-dimensional inviscid Burgers' equation.}
    \label{fig:apostbe2}
\end{figure}

\begin{figure}
    \centering
    \begin{minipage}{\textwidth}
        \centering
        \begin{subfigure}[b]{0.4\textwidth}
            \centering
            \begin{tikzpicture}
\begin{axis}[
axis equal image,
xmin=-0.9416904449462891,
xmax=0.9417643547058105,
ymin=-0.8028290867805481,
ymax=0.987214982509613,
xtick={-.9,0,0.9},
ytick={-.8,0,0.9},
width=1.0\textwidth,]
\addplot []
graphics [xmin=-0.9416904449462891,xmax=0.9417643547058105,ymin=-0.8028290867805481,ymax=0.987214982509613] {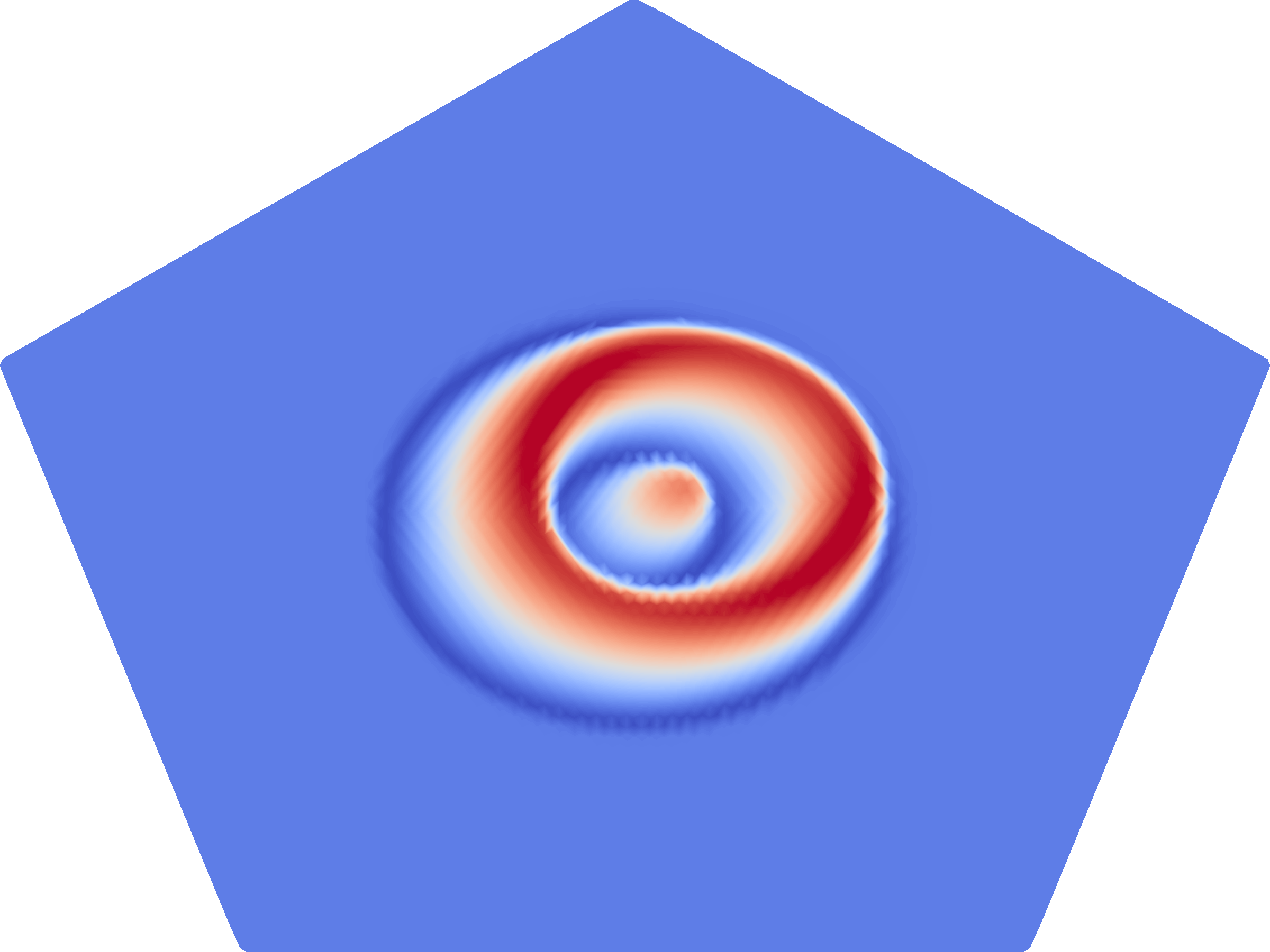};

\end{axis}
\end{tikzpicture}
            \\
            \hspace{2em}\colorbarMatplotlibCoolwarm{0}{0.25}{0.50}{0.75}{1.0}{5cm}{0.3cm} 
            \captionsetup{labelformat=empty}
            \caption{\hspace{0.6cm}(a) Solution at time $t = 0.2\,s$}
            \label{fig:apostbe31}
        \end{subfigure}
        \hspace{0.01\textwidth}
        \begin{subfigure}[b]{0.4\textwidth}
            \centering
            \begin{tikzpicture}
\begin{axis}[
axis equal image,
xmin=-0.9416904449462891,
xmax=0.9417643547058105,
ymin=-0.8028290867805481,
ymax=0.987214982509613,
xtick={-.9, 0, 0.9},
ytick={-.8, 0, 0.9},
width=1.0\textwidth]
\addplot []
graphics [xmin=-0.9416904449462891,xmax=0.9417643547058105,ymin=-0.8028290867805481,ymax=0.987214982509613] {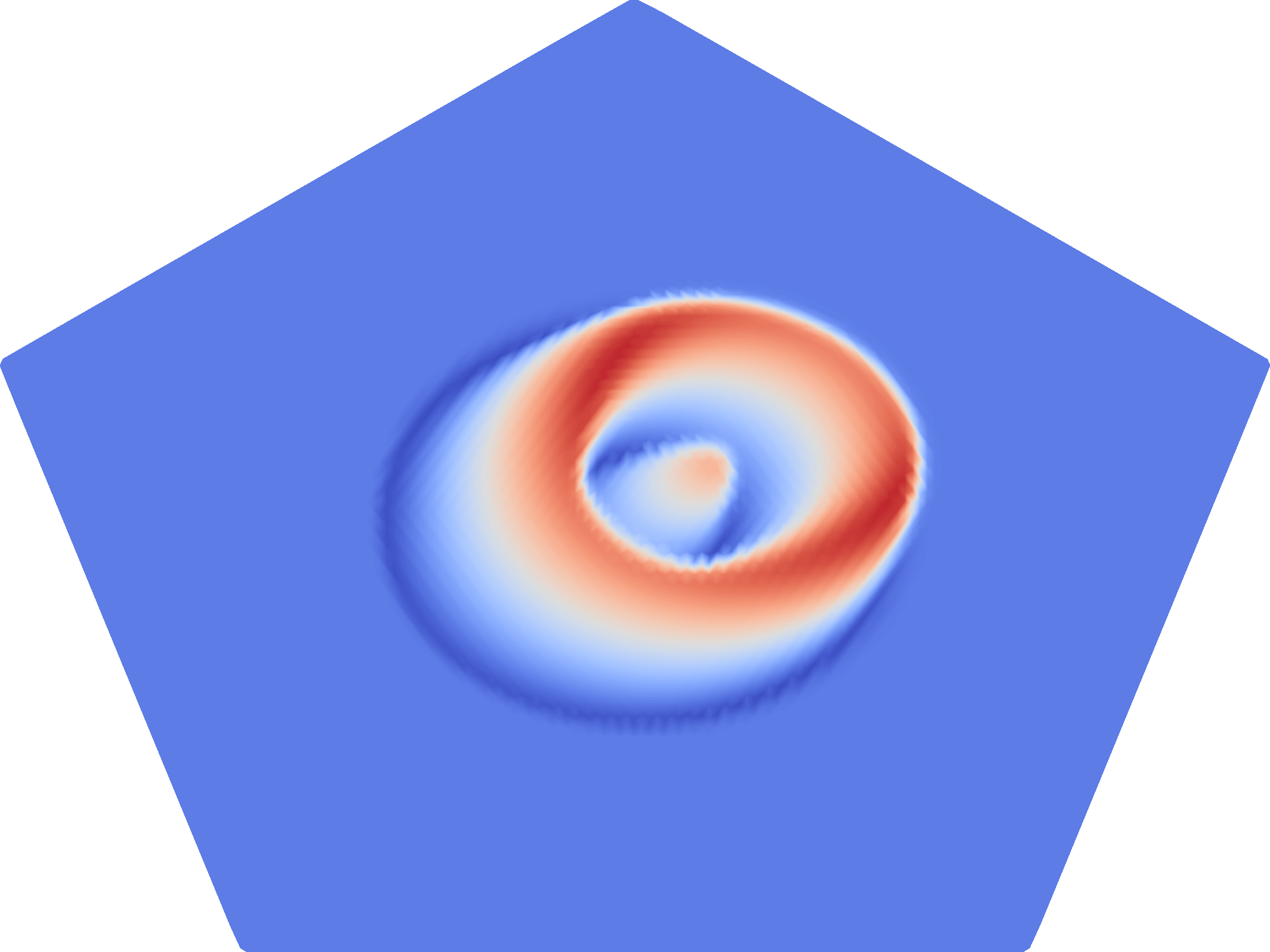};

\end{axis}
\end{tikzpicture}
            \\
            \hspace{2em}\colorbarMatplotlibCoolwarm{0}{0.25}{0.50}{0.75}{1.0}{5cm}{0.3cm} 
            \captionsetup{labelformat=empty}
            \caption{\hspace{0.7cm}(b) Solution at time $t = 0.4\,s$}
            \label{fig:apostbe32}
        \end{subfigure}
    \end{minipage}

    \vspace{1cm}

    \begin{minipage}{\textwidth}
        \centering
        \begin{subfigure}[b]{0.4\textwidth}
            \centering
            \begin{tikzpicture}
\begin{axis}[
axis equal image,
xmin=-0.9416904449462891,
xmax=0.9417643547058105,
ymin=-0.8028290867805481,
ymax=0.987214982509613,
xtick={-.9, 0,0.9},
ytick={-.8, 0,0.9},
width=1.0\textwidth]
\addplot []
graphics [xmin=-0.9416904449462891,xmax=0.9417643547058105,ymin=-0.8028290867805481,ymax=0.987214982509613] {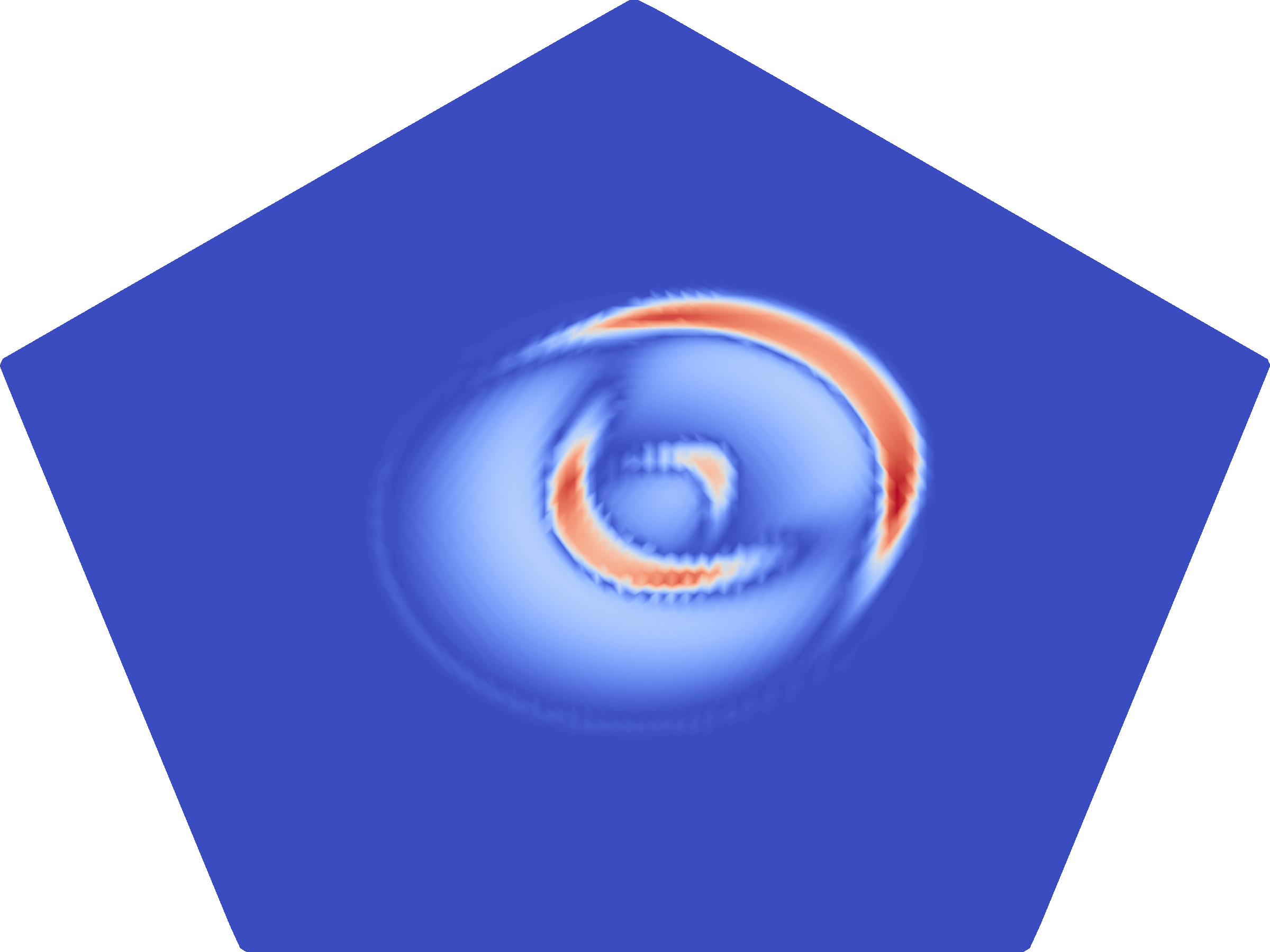};

\end{axis}
\end{tikzpicture}
            \\
            \hspace{2em}\colorbarMatplotlibCoolwarm{0}{0.23}{0.46}{0.69}{0.92}{5cm}{0.3cm} 
            \captionsetup{labelformat=empty}
            \caption{\hspace{0.85cm}(c) Absolute error at time $t = 0.2\,s$}
            \label{fig:apostbe33}
        \end{subfigure}
        \hspace{0.01\textwidth}
        \begin{subfigure}[b]{0.4\textwidth}
            \centering
            \begin{tikzpicture}
\begin{axis}[
axis equal image,
xmin=-0.9416904449462891,
xmax=0.9417643547058105,
ymin=-0.8028290867805481,
ymax=0.987214982509613,
xtick={-.9, 0, 0.9},
ytick={-.8, 0, 0.9},
width=1.0\textwidth]
\addplot []
graphics [xmin=-0.9416904449462891,xmax=0.9417643547058105,ymin=-0.8028290867805481,ymax=0.987214982509613] {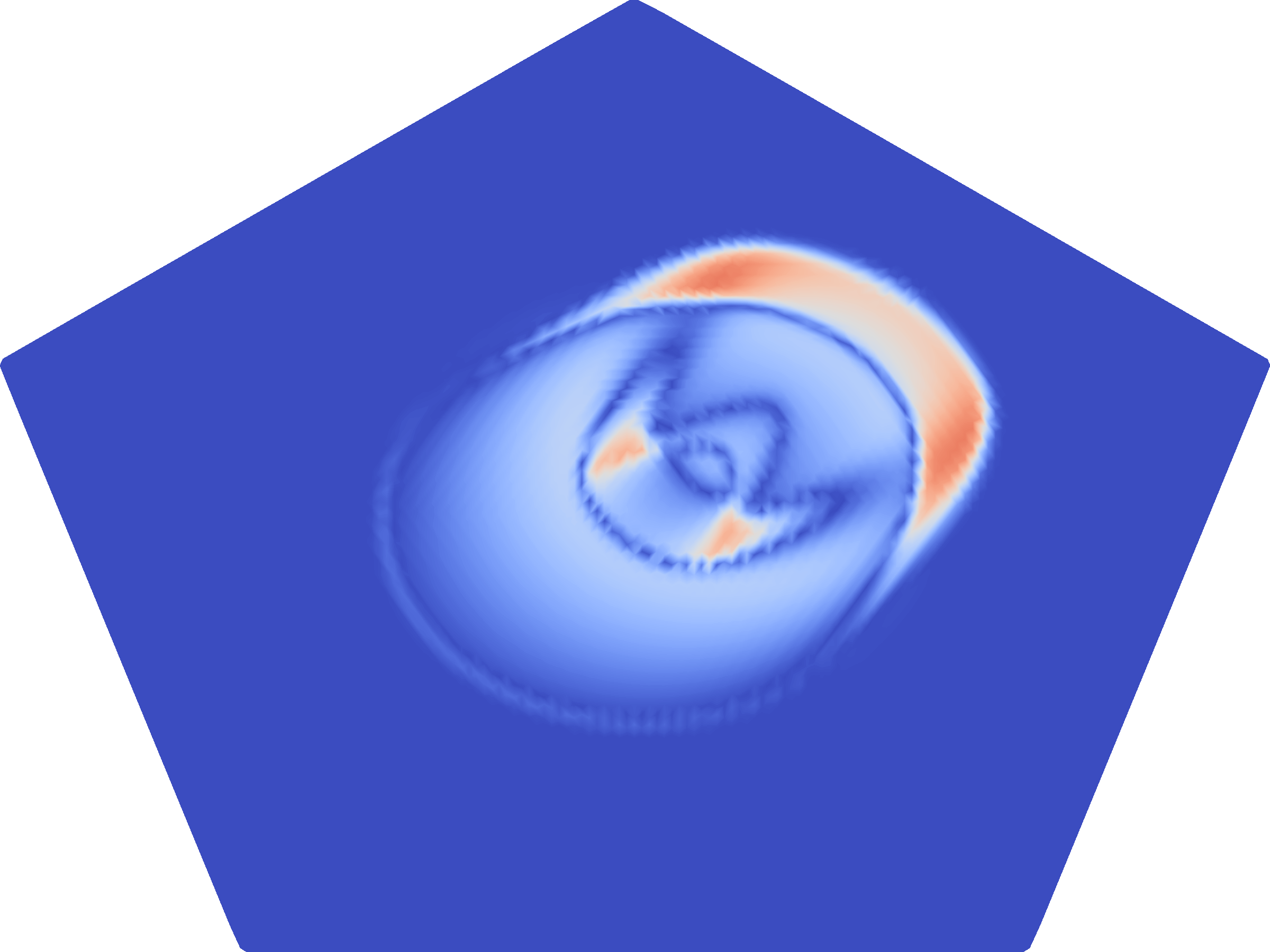};

\end{axis}
\end{tikzpicture}
            \\
            \hspace{2em}\colorbarMatplotlibCoolwarm{0}{0.23}{0.46}{0.69}{0.92}{5cm}{0.3cm} 
            \captionsetup{labelformat=empty}
            \caption{\hspace{0.59cm}(d) Absolute error at time $t = 0.4\,s$}
            \label{fig:apostbe34}
        \end{subfigure}
    \end{minipage}

    \caption{FVM-based numerical solutions obtained using VNN as numerical flux at different time instances and corresponding absolute error field (relative to FVM simulation with Godunov numerical flux) for two-dimensional inviscid Burgers' equation.}
    \label{fig:apostbe3}
\end{figure}

For an a posteriori test of our SM-based fluxes, we consider a two-dimensional regular pentagon as our spatial domain and run FVM simulations for the following initial condition
\[
u(x,0) =
\begin{cases} 
\left(\cos(8 \pi r) + 1\right) \dfrac{e^r}{1 + e^r}, & r < r_\text{max}, \\[0.5em]
\left(\cos(8 \pi r_\text{max}) + 1\right) \dfrac{e^{r_\text{max}}}{1 + e^{r_\text{max}}}, & r \geq r_\text{max},
\end{cases}
\]
where $r=\sqrt{x_1^2+x_2^2}$ is radial distance from the origin, and $r_\text{max} = 0.4$ specifies the cutoff radius. The simulation is carried out over the time interval $\mathcal{T} = (0,0.4)$ using the SMs trained for the one-dimensional Burgers' equation in Section \ref{ssec:burg1d}. The mesh of the pentagonal domain consists of 11049 triangular elements, and the solutions obtained using different SM-based flux functions at $t=0.2\,s$ and $t=0.4\,s$ are presented in Figures \ref{fig:apostbe2} and \ref{fig:apostbe3} along with the corresponding absolute errors with respect to the solution obtained using Godunov flux. Comparing the absolute error plots in Figures \ref{fig:apostbe2} and \ref{fig:apostbe3} reveals that the numerical solutions obtained using the BFNN-based flux function are far more accurate than solutions obtained using VNN. 

\subsection{Two-dimensional shallow water equations} \label{ssec:swe2d}
As our final example, we consider the two-dimensional shallow water equations confined to the domain $\Omega_x \subset \mathbb{R}^2$, which can be expressed as 
\begin{equation}
   \check{u} = \begin{bmatrix} u \\[0.3em]
     v \end{bmatrix}, \quad U = \begin{bmatrix} h \\[0.3em]
     h\check{u} \end{bmatrix}, \quad F_I(U) = \begin{bmatrix} h\check{u}^T \\[0.3em]
     h\check{u}\check{u}^T + \dfrac{h^2}{2}I \end{bmatrix}, \quad F_{V}(U,\nabla U) = 0, \quad S(U,\nabla U) = 0
    \label{eqn:swe2d}
\end{equation}
\noindent where $h(x,t) : \Omega_x \to \mathbb{R}_{>0}$ is the height of the free surface and $\check{u}(x,t): \Omega_x \to \mathbb{R}^2$ is the fluid velocity with $x:= (x_1,x_2) \in \Omega_x$. Again, we reuse the SM developed for one-dimensional SWEs. 

For an a posteriori test of our SM-based fluxes, we conduct FVM simulations on a two-dimensional square spatial domain $\Omega_x = (0,10)^2$ using the following initial condition
\[
\check{u}(x,0) = \begin{bmatrix} 0.0 \\
     0.0 \end{bmatrix} \quad \quad 
h(x,0) =
\begin{cases} 
3.0 & \text{if } (x_1 - 5)^2 + (x_2 - 5)^2 \leq 6.25, \\
0.25 & \text{otherwise}.
\end{cases} 
\]

\begin{figure}
    \centering
    \begin{minipage}{\textwidth}
        \centering
        \begin{subfigure}[b]{0.4\textwidth}
            \centering
            \begin{tikzpicture}
\begin{axis}[
axis equal image,
xmin=0.06329113990068436,
xmax=9.936708450317383,
ymin=0.06329113990068436,
ymax=9.936708450317383,
width=0.55\textwidth,
xtick={1,5, 9},
ytick={1,5,9},
width=1.\textwidth]
\addplot []
graphics [xmin=0.06329113990068436,xmax=9.936708450317383,ymin=0.06329113990068436,ymax=9.936708450317383] {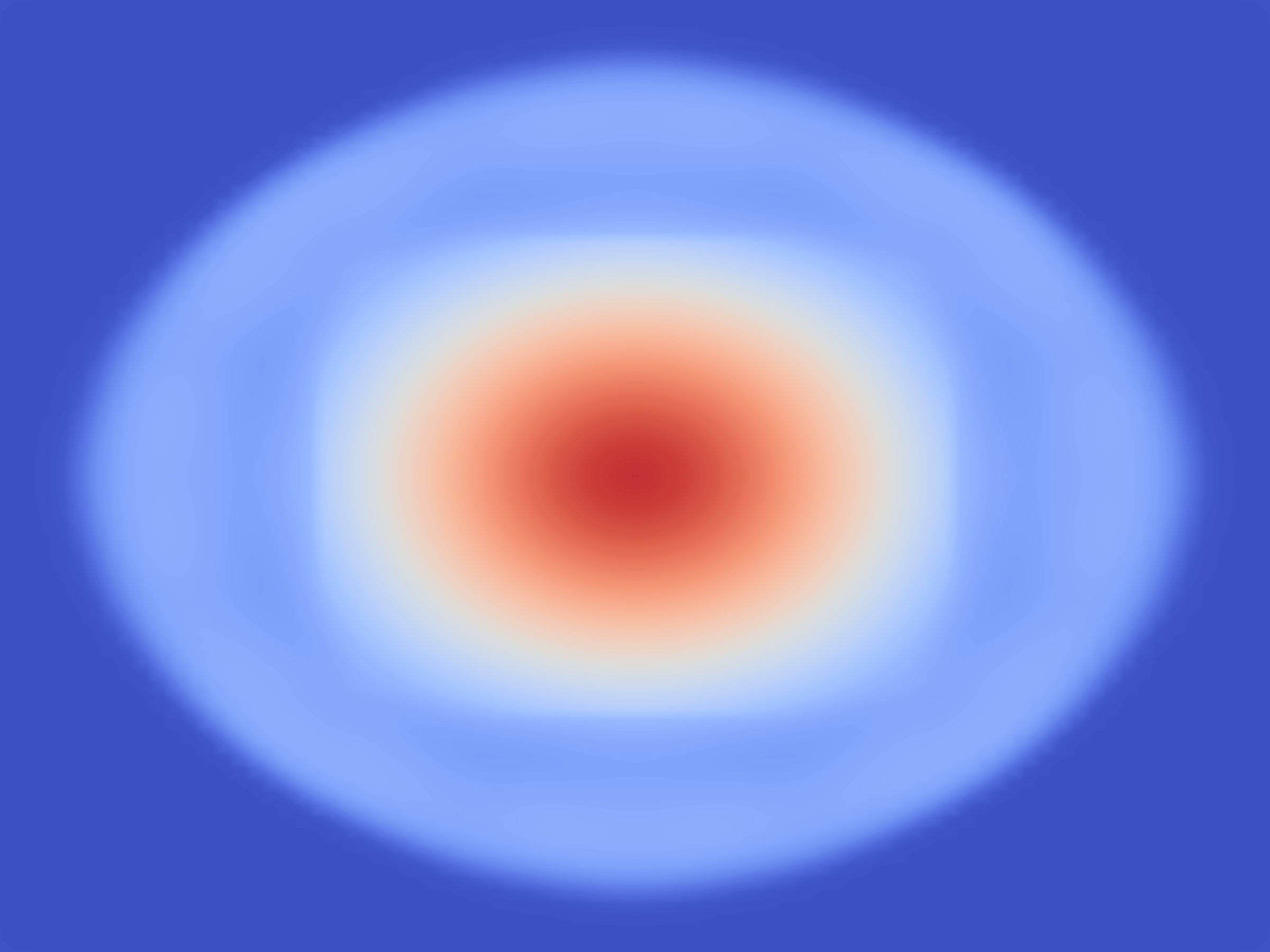};

\end{axis}
\end{tikzpicture}
            \\
            \hspace{1.5em}\colorbarMatplotlibCoolwarm{0.2}{0.8}{1.6}{2.4}{3.0}{5cm}{0.3cm} 
            \captionsetup{labelformat=empty}
            \caption{\hspace{0.6cm}(a) Solution at time $t = 1\,s$}
            \label{fig:apostbe41}
        \end{subfigure}
        \hspace{0.01\textwidth}
        \begin{subfigure}[b]{0.4\textwidth}
            \centering
            \begin{tikzpicture}
\begin{axis}[
axis equal image,
xmin=0.06329113990068436,
xmax=9.936708450317383,
ymin=0.06329113990068436,
ymax=9.936708450317383,
width=0.55\textwidth,
xtick={1,5,9},
ytick={1,5,9},
width=1.\textwidth]
\addplot []
graphics [xmin=0.06329113990068436,xmax=9.936708450317383,ymin=0.06329113990068436,ymax=9.936708450317383] {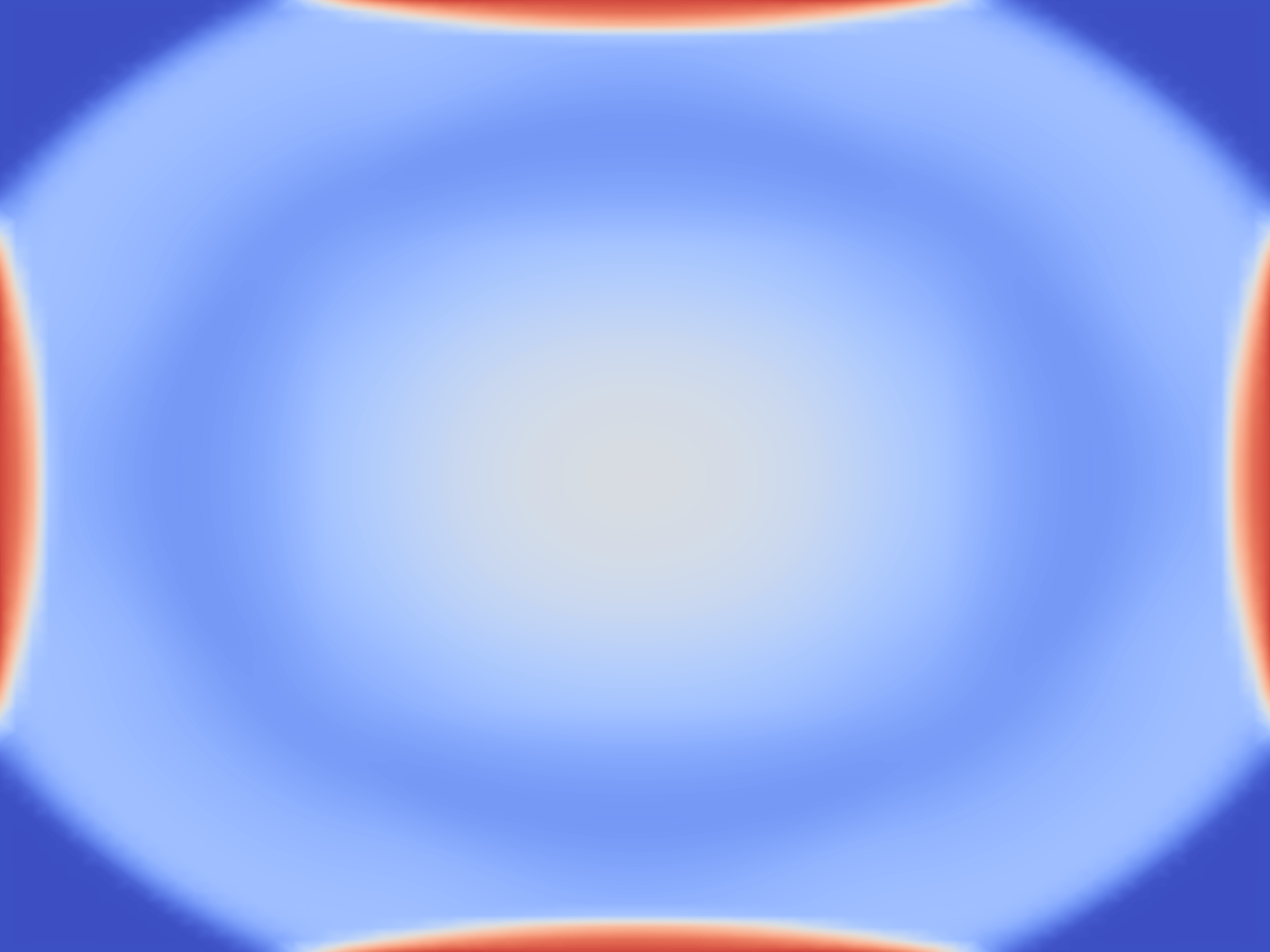};

\end{axis}
\end{tikzpicture}
            \\
            \hspace{1.5em}\colorbarMatplotlibCoolwarm{0.2}{0.7}{1.2}{1.7}{2.2}{5cm}{0.3cm} 
            \captionsetup{labelformat=empty}
            \caption{\hspace{0.7cm}(b) Solution at time $t = 2\,s$}
            \label{fig:apostbe42}
        \end{subfigure}
    \end{minipage}
    
    \vspace{1cm}
    
    \begin{minipage}{\textwidth}
        \centering
        \begin{subfigure}[b]{0.4\textwidth}
            \centering
            \begin{tikzpicture}
\begin{axis}[
axis equal image,
xmin=0.06329113990068436,
xmax=9.936708450317383,
ymin=0.06329113990068436,
ymax=9.936708450317383,
xtick={1,5,9},
ytick={1,5,9},
width=1.0\textwidth]
\addplot []
graphics [xmin=0.06329113990068436,xmax=9.936708450317383,ymin=0.06329113990068436,ymax=9.936708450317383] {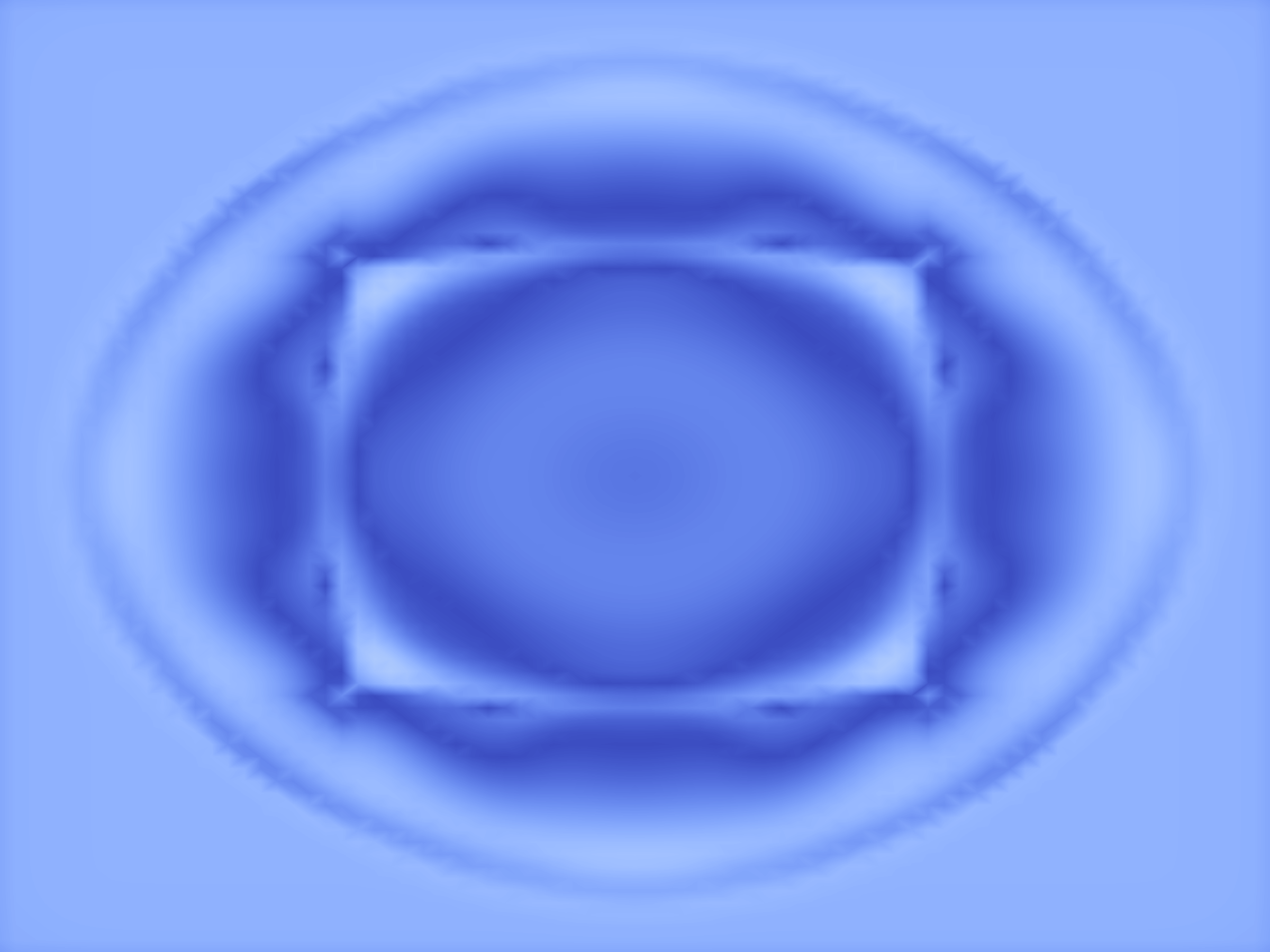};

\end{axis}
\end{tikzpicture}
            \\
            \hspace{1.5em}\colorbarMatplotlibCoolwarm{0}{0.022}{0.046}{0.068}{0.09}{5cm}{0.3cm} 
            \captionsetup{labelformat=empty}
            \caption{\hspace{0.85cm}(c) Absolute error at time $t = 1\,s$}
            \label{fig:apostbe43}
        \end{subfigure}
        \hspace{0.01\textwidth}
        \begin{subfigure}[b]{0.4\textwidth}
            \centering
            \begin{tikzpicture}
\begin{axis}[
axis equal image,
xmin=0.06329113990068436,
xmax=9.936708450317383,
ymin=0.06329113990068436,
ymax=9.936708450317383,
xtick={1,5,9},
ytick={1,5,9},
width=1.0\textwidth]
\addplot []
graphics [xmin=0.06329113990068436,xmax=9.936708450317383,ymin=0.06329113990068436,ymax=9.936708450317383] {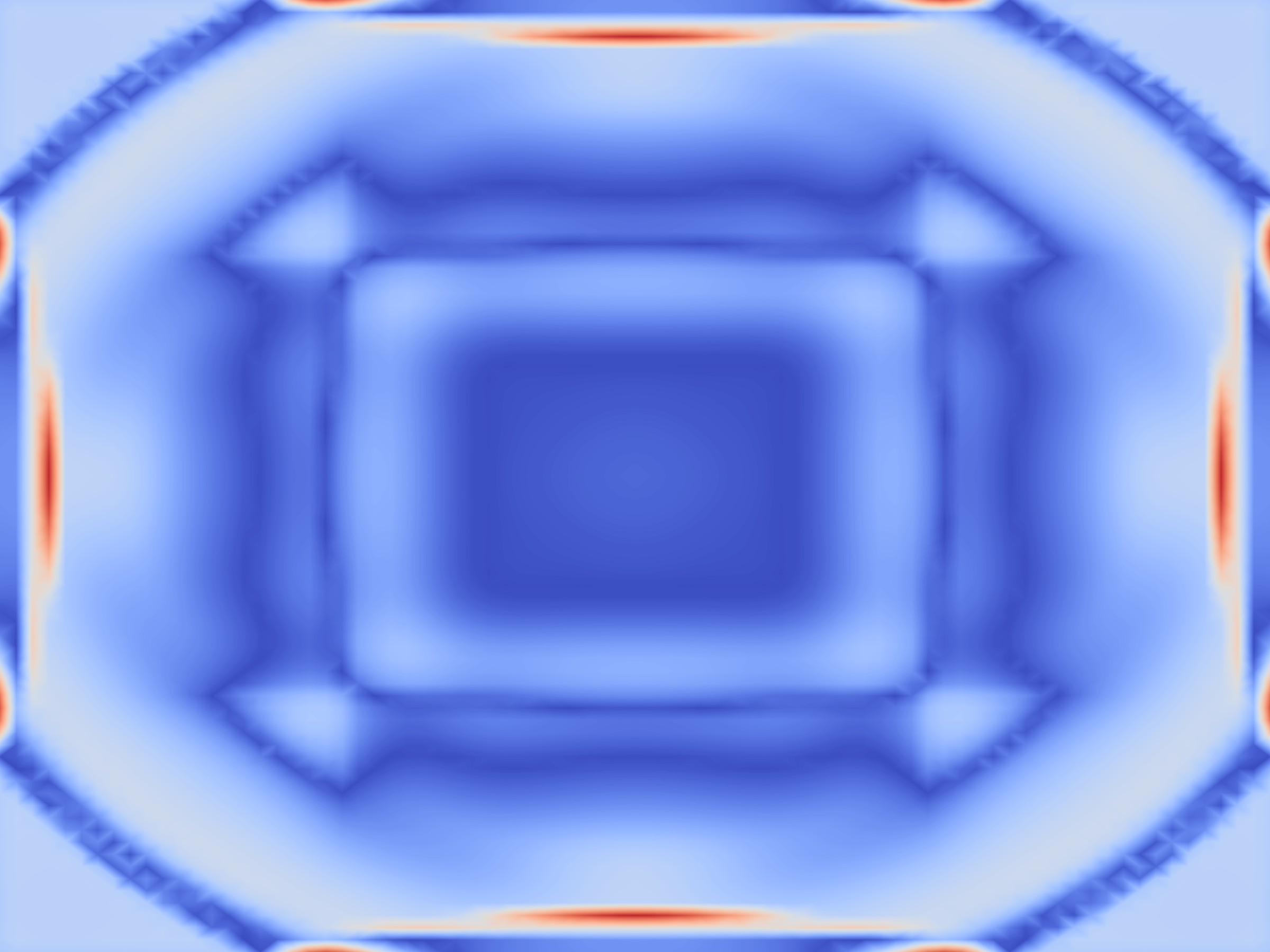};

\end{axis}
\end{tikzpicture}
            \\
            \hspace{1.5em}\colorbarMatplotlibCoolwarm{0}{0.022}{0.046}{0.068}{0.09}{5cm}{0.3cm} 
            \captionsetup{labelformat=empty}
            \caption{\hspace{0.59cm}(d) Absolute error at time $t = 2\,s$}
            \label{fig:apostbe44}
        \end{subfigure}
    \end{minipage}

    \caption{FVM-based numerical solutions (height) obtained using BFNN as numerical flux at different time instances and corresponding absolute error field (relative to FVM simulation with Godunov numerical flux) for the two-dimensional shallow water equations.}
    \label{fig:apostbe4}
\end{figure}

\begin{figure}
    \centering
    \begin{minipage}{\textwidth}
        \centering
        \begin{subfigure}[b]{0.4\textwidth}
            \centering
            \begin{tikzpicture}
\begin{axis}[
axis equal image,
xmin=0.06329113990068436,
xmax=9.936708450317383,
ymin=0.06329113990068436,
ymax=9.936708450317383,
xtick={1,5,9},
ytick={1,5,9},
width=1.0\textwidth]
\addplot []
graphics [xmin=0.06329113990068436,xmax=9.936708450317383,ymin=0.06329113990068436,ymax=9.936708450317383] {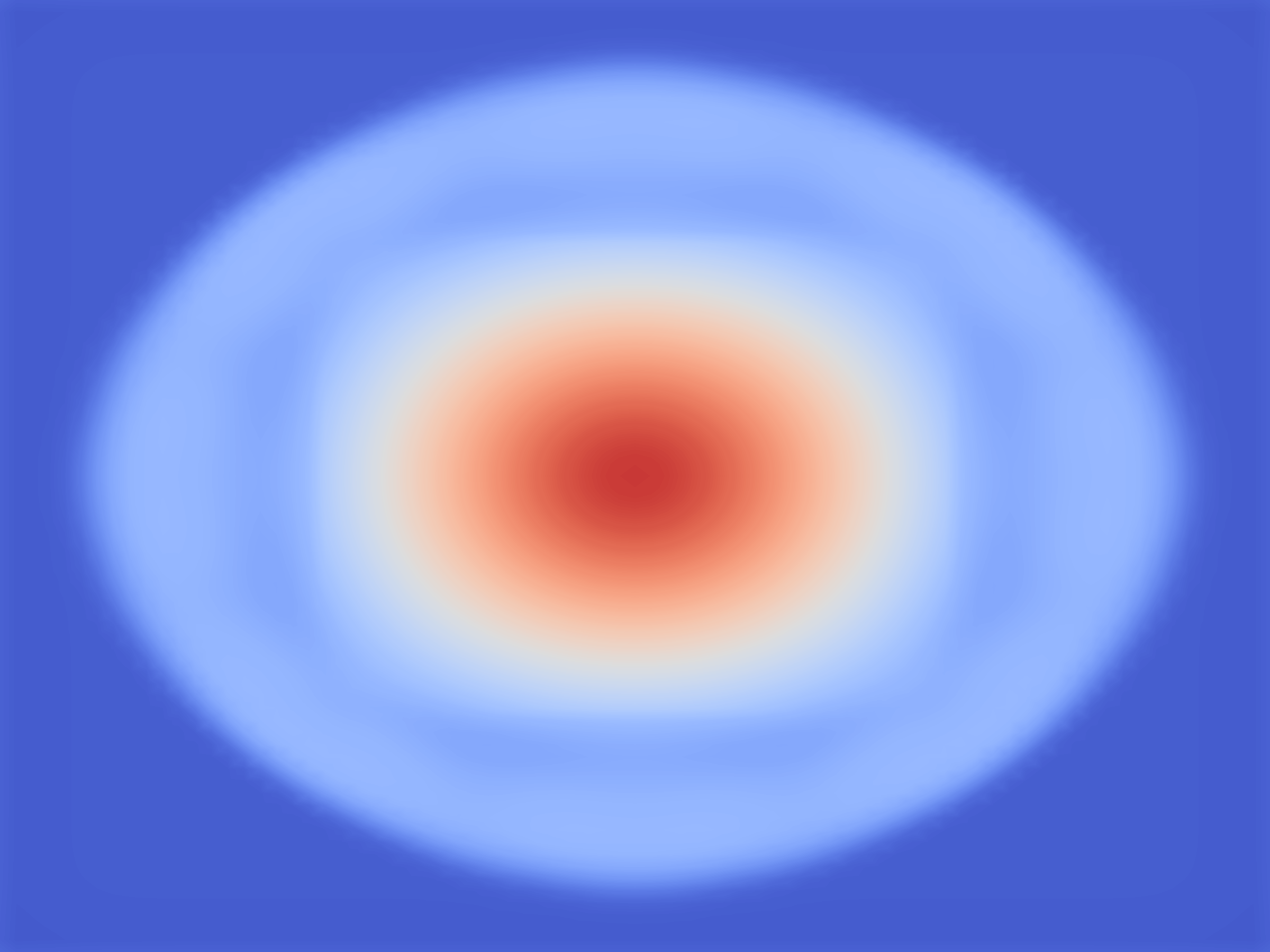};

\end{axis}
\end{tikzpicture}
            \\
            \hspace{1.5em}\colorbarMatplotlibCoolwarm{0.2}{0.8}{1.6}{2.4}{3.0}{5cm}{0.3cm} 
            \captionsetup{labelformat=empty}
            \caption{\hspace{0.6cm}(a) Solution at time $t = 1\,s$}
            \label{fig:apostbe51}
        \end{subfigure}
        \hspace{0.01\textwidth}
        \begin{subfigure}[b]{0.4\textwidth}
            \centering
            \begin{tikzpicture}
\begin{axis}[
axis equal image,
xmin=0.06329113990068436,
xmax=9.936708450317383,
ymin=0.06329113990068436,
ymax=9.936708450317383,
xtick={1,5,9},
ytick={1,5,9},
width=1.0\textwidth]
\addplot []
graphics [xmin=0.06329113990068436,xmax=9.936708450317383,ymin=0.06329113990068436,ymax=9.936708450317383] {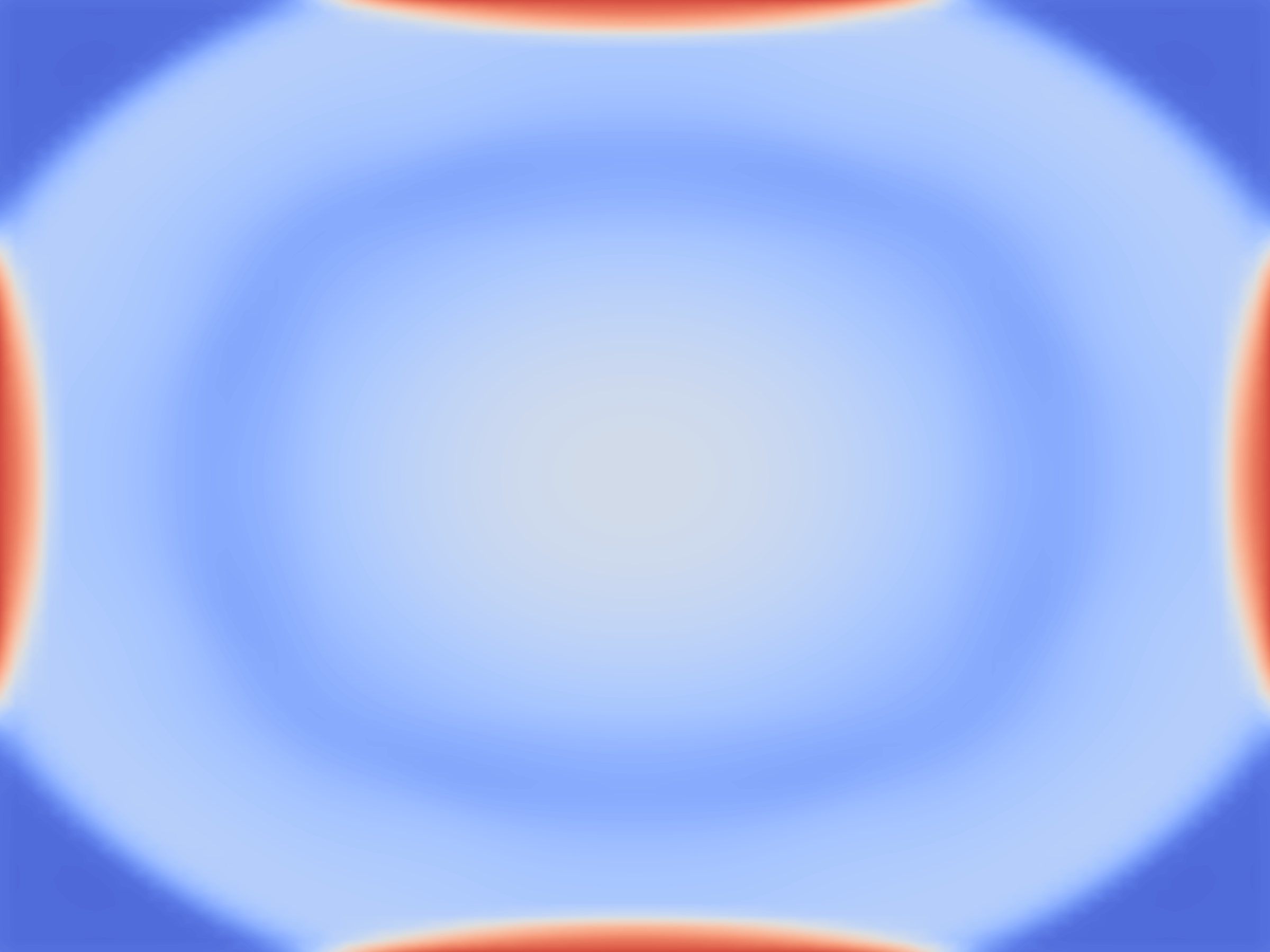};

\end{axis}
\end{tikzpicture}
            \\
            \hspace{1.5em}\colorbarMatplotlibCoolwarm{0.2}{0.7}{1.2}{1.7}{2.2}{5cm}{0.3cm} 
            \captionsetup{labelformat=empty}
            \caption{\hspace{0.7cm}(b) Solution at time $t = 2\,s$}
            \label{fig:apostbe52}
        \end{subfigure}
    \end{minipage}
    
    \vspace{1cm}
    
    \begin{minipage}{\textwidth}
        \centering
        \begin{subfigure}[b]{0.4\textwidth}
            \centering
            \begin{tikzpicture}
\begin{axis}[
axis equal image,
xmin=0.06329113990068436,
xmax=9.936708450317383,
ymin=0.06329113990068436,
ymax=9.936708450317383,
xtick={1,5,9},
ytick={1,5,9},
width=1.0\textwidth]
\addplot []
graphics [xmin=0.06329113990068436,xmax=9.936708450317383,ymin=0.06329113990068436,ymax=9.936708450317383] {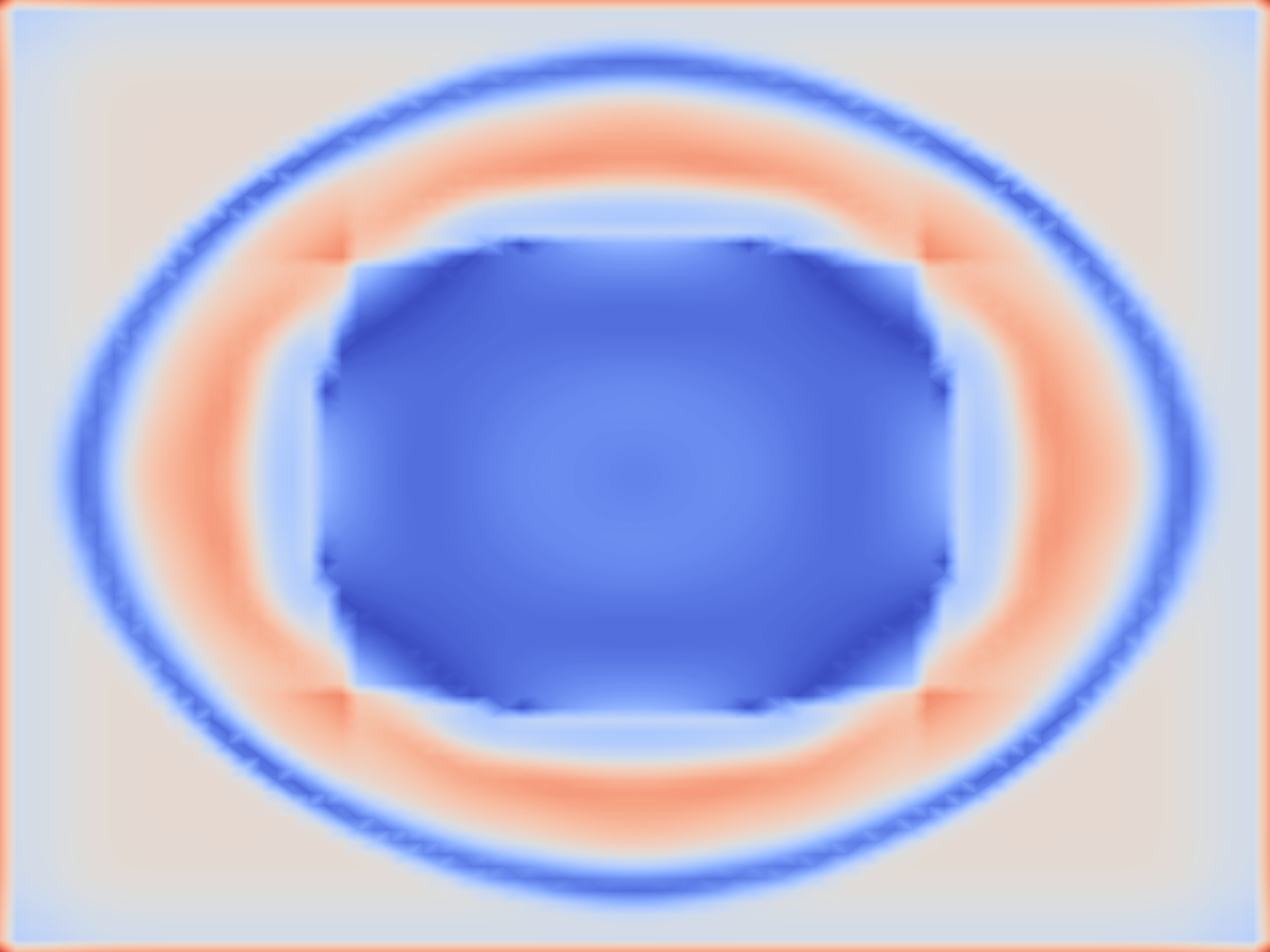};

\end{axis}
\end{tikzpicture}
            \\
            \hspace{1.5em}\colorbarMatplotlibCoolwarm{0.00}{0.03}{0.06}{0.09}{0.12}{5cm}{0.3cm} 
            \captionsetup{labelformat=empty}
            \caption{\hspace{0.85cm}(c) Absolute error at time $t = 1\,s$}
            \label{fig:apostbe53}
        \end{subfigure}
        \hspace{0.01\textwidth}
        \begin{subfigure}[b]{0.4\textwidth}
            \centering
            \begin{tikzpicture}
\begin{axis}[
axis equal image,
xmin=0.06329113990068436,
xmax=9.936708450317383,
ymin=0.06329113990068436,
ymax=9.936708450317383,
xtick={1,5,9},
ytick={1,5,9},
width=1.0\textwidth]
\addplot []
graphics [xmin=0.06329113990068436,xmax=9.936708450317383,ymin=0.06329113990068436,ymax=9.936708450317383] {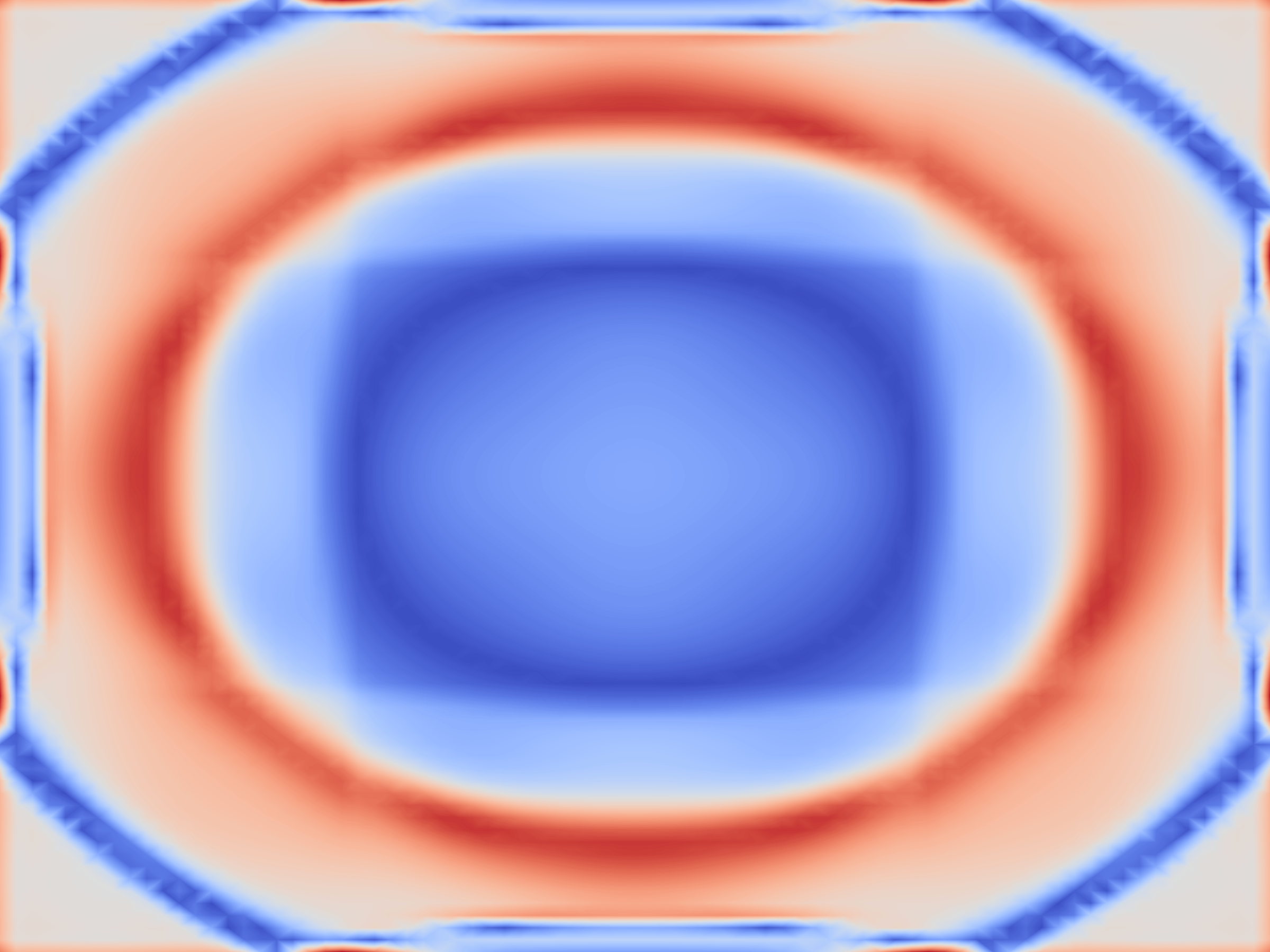};

\end{axis}
\end{tikzpicture}
            \\
            \hspace{1.5em}\colorbarMatplotlibCoolwarm{0.00}{0.04}{0.09}{0.13}{0.18}{5cm}{0.3cm} 
            \captionsetup{labelformat=empty}
            \caption{\hspace{0.59cm}(d) Absolute error at time $t = 2\,s$}
            \label{fig:apostbe54}
        \end{subfigure}
    \end{minipage}

	\caption{FVM-based numerical solutions (height) obtained using VNN as numerical flux at different time instances and corresponding absolute error field (relative to FVM simulation with Godunov numerical flux) for the two-dimensional shallow water equations.}
    \label{fig:apostbe5}
\end{figure}

 SMs trained for one-dimensional SWEs are used as numerical flux functions and simulations are performed over the time interval $\mathcal{T} = (0, 2.0)$. The computational mesh for the square domain consists of 6400 quadrilateral elements. The solutions computed using the two SM-based flux functions at $t = 1.0\,\text{s}$ and $t = 2.0\,\text{s}$ are presented in Figures \ref{fig:apostbe4} and \ref{fig:apostbe5}. These figures also contain the corresponding absolute error plots, where the error is evaluated relative to the reference solution obtained with the Godunov flux. Similar to the results obtained for previous examples, we find that the solutions obtained using BFNN as a numerical flux function are more accurate as compared to VNN.

\section{Conclusion}
\label{sec:conclude}
In this work, we propose Godunov flux surrogates using NN-based approximations that map interior and exterior conservative state variables to the corresponding flux. The objective is to retain the robustness and accuracy inherent to the Godunov flux. We developed two distinct NN-based surrogate models: the first utilized a vanilla FCNN trained via supervised learning, and the second leveraged a bi-fidelity NN-based model. In the bi-fidelity approach, an FCNN was trained to correct an approximate numerical flux function, with the residual between the Godunov flux and the approximate flux serving as the target for supervised learning. We constructed the bi-fidelity model using the Roe numerical flux function as the LF solver. We demonstrated the performance of the proposed surrogate models through applications to both one-dimensional and two-dimensional PDEs. For the one-dimensional case, we conducted a priori tests on the inviscid Burgers' equation and shallow water equations, as the surrogate models were specifically trained on one-dimensional Godunov flux data. In contrast, a posteriori tests were performed on all example PDEs, including two-dimensional cases. Our results revealed that the BFNN-based SM exhibited superior accuracy and generalizability compared to the VNN-based SM. Moreover, we also presented a priori and a posteriori test cases, where the proposed SMs were more accurate and robust than approximate numerical flux functions (e.g., Roe, Roe with Harten fix, and HLL). Lastly, we exploited symmetries inherent in conservation laws to adapt the trained one-dimensional surrogate models for use as flux functions in two-dimensional FVM simulations. However, in this study, the training datasets were generated by sampling the input states $U^+$ and $U^-$ from uniform distributions. While this choice ensures broad coverage of the input space, future work could explore more targeted sampling strategies, such as selecting the input states from distributions that lead to specific wave structures in the Riemann solution, to potentially further enhance model performance and robustness in particular regimes of interest.

\section*{Acknowledgments}
This work is supported by
AFOSR award number FA9550-20-1-0236 and
NSF award number CBET-2338843.
The content of this publication does not necessarily reflect the position
or policy of any of these supporters, and no official endorsement should
be inferred.

\bibliographystyle{plain}
\bibliography{biblio}

\end{document}